\documentclass[review,12pt]{elsarticle}

\usepackage{amssymb}

\usepackage{lineno}

\usepackage{natbib}
\usepackage{geometry}
\usepackage{amsmath,amsfonts,amsthm}
\usepackage[utf8x]{inputenc} 
\usepackage[english]{babel} 
\usepackage[T1]{fontenc} 
\usepackage{latexsym}
\usepackage{psfrag}
\usepackage{epsfig}
\usepackage{graphicx}
\usepackage{exscale}
\usepackage{url} 
\usepackage{textcomp}
\usepackage{float}
\usepackage[bookmarks,colorlinks]{hyperref}
\usepackage{amscd}
\usepackage{indentfirst}
\usepackage{graphicx}
\usepackage{bm}
\usepackage[]{subfig}
\usepackage{fancyhdr}
\usepackage{booktabs}
\usepackage{caption}
\usepackage{physics}
\usepackage{color}
\definecolor{hillencolor}{rgb}{0.66, 0.030, 0.30}

\newtheorem{Theorem}{Theorem}[section]

\newtheorem{Proposition}{Proposition}[section]

\newenvironment{Proof}{\par\noindent{\bf Proof.}}{\hfill$\square$}

\begin{document}
\begin{frontmatter}



\title{Weakly nonlinear analysis of a two-species non-local advection-diffusion system}

\author{Valeria Giunta\corref{cor1}\fnref{label1}}
\cortext[cor1]{Corresponding author}
\author[label2]{Thomas Hillen}
\author[label3]{Mark A. Lewis}
\author[label1]{Jonathan R. Potts}
 \affiliation[label1]{School of Mathematics and Statistics, University of Sheffield, Hicks Building, Hounsfield Road,
Sheffield S3 7RH, UK}
 \affiliation[label2]{Department of Mathematical and Statistical Sciences, University of Alberta, Edmonton,
AB T6G 2G1, Canada}
 \affiliation[label3]{Department of Mathematics and Statistics and Department of Biology, University of Victoria,
PO Box 1700 Station CSC, Victoria, BC, Canada}



\begin{abstract}
Nonlocal interactions are ubiquitous in nature and play a central role in many biological systems. In this paper, we perform a bifurcation analysis of a widely-applicable advection-diffusion model with nonlocal advection terms describing the species movements generated by inter-species  interactions. We use linear analysis to assess the stability of the constant steady state, then weakly nonlinear analysis to recover the shape and stability of non-homogeneous solutions. Since the system arises from a conservation law, the resulting amplitude equations consist of a Ginzburg-Landau equation coupled with an equation for the zero mode. In particular, this means that supercritical branches from the Ginzburg-Landau equation need not be stable.  Indeed, we find that, depending on the parameters, bifurcations can be subcritical (always unstable), stable supercritical, or unstable supercritical.  We show numerically that, when small amplitude patterns are unstable, the system exhibits large amplitude patterns and hysteresis, even in supercritical regimes.  Finally, we construct bifurcation diagrams by combining our analysis with a previous study of the minimisers of the associated energy functional. Through this approach we reveal parameter regions in which stable small amplitude patterns coexist with strongly modulated solutions.

\end{abstract}




\end{frontmatter}


\section{Introduction}

    Spontaneous pattern formation occurs throughout nature \cite{hoyle06}, with examples ranging from animal coat patterns \cite{Turing1952} to territory formation \cite{potts2016territorial}, cell sorting \cite{Stevens} and swarm aggregation \cite{TBL06}.
    Therefore uncovering and analysing the mechanisms behind pattern formation 
    is a central challenge in the life sciences where applied mathematics can play a role.  Typically, research into pattern formation proceeds first by assessing which parameters may cause patterns to emerge spontaneously from a homogeneous steady state, using linear pattern formation analysis, sometimes called `Turing pattern analysis' \cite{Turing1952}.  This determines whether patterns may emerge at short times from arbitrarily small perturbations.  However, it is also important biologically to show whether these patterns are stable. One approach to pattern stability is via weakly nonlinear analysis: a stable supercritical bifurcation branch suggest that asymptotic patterns will emerge continuously as the bifurcation parameter is changed, whereas an unstable subcritical branch suggests that large amplitude asymptotic patterns may  appear abruptly as the bifurcation point is crossed, their amplitude being a discontinuous function of the bifurcation parameter.  This discontinuity in amplitude with respect to parameter change indicates that a biological system might suddenly change its behaviour in a dramatic fashion with only a small change in the underlying mechanisms.

    Many biological mechanisms generate attractive or repulsive forces governing phenomena such as chemotaxis (\cite{eisenbach2004chemotaxis,usersguide}), bacterial orientation (\cite{basaran2022large}), swarms of animals (\cite{puckett2014searching}), and motion of human crowds (\cite{goldstone2009collective}). These mechanisms are driven by electrical, chemical or social interactions. These interactions arise from individual organisms collecting information from their environment, such as the presence of other individuals, food or chemicals. After gathering information, individuals move towards regions that contain important components for survival or move away from less favourable areas, thus creating spatially inhomogeneous distributions of individuals, which may have a certain degree of regularity in space and/or time (e.g. \cite{TBL06, PL19}). This process of acquiring information from the environment is generally nonlocal, as motile organisms are usually able to inspect a portion of their environment, either by prolonging their protrusions, as in the case of cells \cite{buttenschon2021non}, or by using their sight, hearing or smell, as with animals \cite{potts2014animal}. 
    
    In recent years there has been an increasing interest in the mathematical modelling of nonlocal advection as a movement model with nonlocal information \cite{BLP02,TBL06,Stevens,CKY19,buttenschon2021non}. Recently, the following class of nonlocal advection-diffusion equations was proposed as a general model of interacting  populations \cite{PL19} 
    \begin{linenomath*}\begin{equation}\label{eq:system}
		\frac{\partial u_i}{\partial t}=D_i \Delta u_i+\nabla \cdot \left(u_i \sum_{j=1}^{N}\gamma_{i j} \nabla (K \ast u_j) \right), \, i=1,\dots, N.
    \end{equation}\end{linenomath*}
    Here, $ u_i(x,t) $ denotes the density of population $ i $ at position $ x $ and time $ t $, for $ i \in \{1, \dots, N\} $ and   $ D_i >0 $ is the diffusion rate of $ u_i $. Individuals can detect the presence of other individuals, whether conspecifics or not, over a spatial neighborhood described by spatial averaging kernel $K$, which is a symmetric, non-negative function modelling the sensing range. The term $ K \ast u_j $ denotes the convolution between $K$ and $u_j$ and describes the nonlocal interactions of $u_i$ with $u_j$. The parameters $ \gamma_{i j} $ are the inter/intra-species interaction parameters, giving the density-dependent rate at which species $i$ advects towards (if $ \gamma_{i j} <0$), or away from (if $ \gamma_{i j} >0$), species $j$.

    Model \eqref{eq:system} implicitly focuses on time scales whereby birth and death processes are negligible. Nonetheless, it has a wide range of possible applications in that it generalizes a variety of existing models describing many different phenomena, such as animal home ranges \cite{BLP02}, territory formation \cite{ellison2020mechanistic, potts2016territorial, rodriguez2020steady}, and cell sorting \cite{Stevens}. 
    On the mathematical side, well-posedness of System \eqref{eq:system} was analyzed in \cite{giunta2022local} and \cite{jungel2022nonlocal}. When the kernel $K$ is sufficiently smooth, \cite{giunta2022local} shows that the system admits classical, positive and global solutions in 1D dimension, and local strong solutions in any higher dimension. When the kernel is non-smooth, in \cite{jungel2022nonlocal} it is proven that System \eqref{eq:system} has weak solutions that exist globally in time.
    
    From the perspective of pattern formation, numerical analysis shows that System \eqref{eq:system} exhibits a great variety of spatio-temporal patterns, depending on the model parameters.  These include segregated and aggregated stationary patterns, periodic time oscillating solutions, and aperiodic spatio-temporal behaviours \cite{PL19}, \cite{giunta2022local}, { \cite{carrillo2020long}. In many cases the system admits an energy functional \cite{GHLP22, carrillo2020long},  which can be used to gain analytic insight into the steady asymptotic patterns that can form from this system.}   Although \cite{GHLP22} focused on the $N=2$ case, the methods are more generally applicable in principle.

    Here, we  perform a bifurcation analysis of one of the cases analyzed in \cite{GHLP22}, namely where $N=2$, $\gamma_{ij}=\gamma_{ji}$ and $\gamma_{ii}=0$. For simplicity, we also assume that $D_1=D_2$. 
    We use weakly nonlinear analysis to derive the equations governing the amplitude of the stationary solutions. Through analysis of the amplitude equations, we determine the nature of bifurcations generating branches of non-homogeneous solutions from a homogeneous state, then recover the shape of the non-homogeneous solutions and their stability. We validate our results through numerical analysis, setting $K$ to be the top-hat distribution \cite{GHLP22}.
    Finally, we  combine our results with results of \cite{GHLP22} that were derived from an energy principle, to construct bifurcation diagrams that incorporate all the existing analysis of this system.

    An interesting feature of our analysis is that the equation governing the modulation of small-amplitude patterns 
    is not always the real Ginzburg-Landau (GL) equation.  This contrasts with many examples of weakly nonlinear analysis, where the GL equation provides the amplitude of the stationary pattern and its stability: in subcritical regimes, the pattern solution is always unstable; in supercritical regimes, a periodic pattern is stable if its wavenumber lies within the Eckhaus band;
    \cite{TB90,hoyle06,bergeon2008eckhaus,bilotta2018eckhaus,giunta2021pattern,consolo2022eckhaus}.   
    In our case, the real GL equation does not always provide a correct description of the pattern near the onset. This is because our system possesses a conservation law, i.e. mass is conserved for all time. This conservation law gives rise to a large-scale neutral mode (the zero mode) that can affect the stability of the pattern, so must be included into the analysis \cite{CH93, MC00}. Therefore, the resulting amplitude equations will consist of the GL equation coupled to an equation for the large-scale mode. 
    
    In \cite{MC00} the authors used symmetry and scaling arguments to derive the amplitude equations governing systems with a conserved quantity. 
   They proved that there exist stable stationary solutions in the form of strongly modulated patterns (i.e. patterns that consist of multiple Fourier modes), and these exist away from the branch that bifurcates from the constant steady state.  The existence of strongly modulated patterns for System \eqref{eq:system} has also been shown in \cite{GHLP22} by analyzing the minimizers of an energy functional associated with the system. Here we build on this by investigating the existence and stability of small amplitude patterns, and showing that when these solutions are unstable, the system evolves towards either large amplitude or strongly modulated patterns. In addition, our analysis shows that, in some parameter regions, stable small amplitude patterns can coexist with stable strongly modulated solutions.
    
   {A similar two-species aggregation model was studied recently in \cite{Stevens}. Their model differs from our model (\ref{eq:system1}) in regard of the diffusion term. In \cite{Stevens} the terms $D\partial_{xx}u_i $ for $i=1,2$ are replaced by density dependent diffusion terms $D\partial_x(u_i\partial_x (u_1+u_2))$. The pattern forming mechanism is similar to our model, however, the arising aggregations have compact support.}   
    
    This paper is organised as follows. Linear stability analysis is given in Section \ref{sec:lin_analysis} and a weakly nonlinear analysis in Section \ref{sec:wnl_analysis}. In these two sections, the analysis is carried out with a generic kernel, in order to provide some general results that can be used for future works. Section \ref{sec:th} focuses on detailed analysis where $K$ is the top-hat distribution.  We analyse the amplitude equations, recover the bifurcation diagrams and compare analytical results with numerical solutions. We finally combine the analysis performed here with the results obtained in \cite{GHLP22} to recover more exhaustive pictures of the bifurcation diagrams. In Section \ref{sec:discussion}, we outline further extensions of this work and discuss possible applications of our results to natural systems.

    \section{Linear stability analysis}\label{sec:lin_analysis}
     We consider System \eqref{eq:system} with two interacting populations, $u_1$ and $u_2$, that either mutually avoid or attract with the same strength (i.e. $ \gamma_{12}=\gamma_{21} $). We set $ \gamma:= \gamma_{12}=\gamma_{21}$ and fix $ D_1= D_2=:D $, and $ \gamma_{11}=\gamma_{22}=0 $. Therefore, System \eqref{eq:system} reads as
	\begin{linenomath*}\begin{equation}\label{eq:system1}
		\begin{aligned}
			&\partial_t u_1= D \partial_{xx} u_1+ \gamma\partial_{x} \left(u_1  \partial_{x}  (K \ast u_2) \right), \\
			&\partial_t u_2= D \partial_{xx} u_2+ \gamma\partial_{x} \left(u_2 \partial_{x}  (K \ast u_1) \right).
		\end{aligned}
	\end{equation}\end{linenomath*}
    We work on the one dimensional spatial domain $ \Omega=\left[-\frac{l}{2},\frac{l}{2}\right] $ and impose periodic boundary conditions
	\begin{linenomath*}\begin{equation}\label{eq:bc}
		u_i\left(-\frac{l}{2},t\right) = u_i\left(\frac{l}{2},t\right), \quad	\partial_x u_i\left(-\frac{l}{2},t\right) = \partial_x u_i\left(\frac{l}{2},t\right), \quad \text{ for } i \in \{1, 2\} \text{ and } t\geq0.
	\end{equation}\end{linenomath*}
    We consider an even and non-negative kernel $K$ such that
        \begin{linenomath*}\begin{equation}
        \label{eq:kint}
          \int_{-l/2}^{l/2} K(x)dx=1, \text{ and } \, \text{Supp} (K)  = \{x \in \mathbb{R}: K(x) > 0\}= [-\alpha, \alpha]
        \end{equation}\end{linenomath*}
    where the constant $\alpha$ denotes the sensitivity radius. We assume that $\alpha < l/2$. Due to the periodic boundary conditions, we also assume that  $K(x)$ is wrapped around periodically over the domain.

    The periodic boundary conditions (Equation \eqref{eq:bc}) ensure that in System \eqref{eq:system1} the total mass of each population $u_i$ is conserved in time. Indeed the following identities are satisfied
	\begin{linenomath*}\begin{equation}
		\frac{d}{dt} \int_{-l/2}^{l/2} u_i(x,t)  \text{d}x=0, \qquad \text{ for } i=1, 2.
	\end{equation}\end{linenomath*}
	Hence 
	\begin{linenomath*}\begin{equation}\label{eq:int_cond}
		\int_{-l/2}^{l/2} u_i (x,t) \text{d} x= \int_{-l/2}^{l/2} u_i(x,0) \text{d} x=:p_i, \text{ for all }  t \geq0,
	\end{equation}\end{linenomath*}
    where the constant $ p_i $ denotes the size of population $u_i$, for $i=1, 2$. 
    
    Equation \eqref{eq:int_cond} implies that system \eqref{eq:system1} has a unique equilibrium point given by
        \begin{linenomath*}\begin{equation}
             \mathbf{\bar{{u}}}:=(\bar{u}_1, \bar{u}_2)=\left(\frac{p_1}{l},\frac{p_2}{l}\right).
        \end{equation}\end{linenomath*}
\subsection{{Nondimensionalization}}
    We start our analysis by rescaling the original system \eqref{eq:system1} using the following non-dimensional coordinates and variables
	\begin{linenomath*}\begin{equation}\label{eq:nd}
	    \tilde{x}=\frac{x}{\alpha}, \quad \tilde{t}=\frac{D}{\alpha^2}t, \quad \tilde{u}_1= l u_1, \quad \tilde{u}_2=l u_2.
	\end{equation}\end{linenomath*}
    Note that, instead of $\alpha$, one could have rescaled using any other constant that is proportional to the standard deviation of $K(x)$ instead, which may be useful if $K(x)$ does not have compact support, for example.
    
    In the non-dimensional spatial domain, we define the following kernel
        \begin{linenomath*}\begin{equation}\label{eq:ndK}
            \tilde{K}(\tilde{x}) :=  \alpha K(\alpha \tilde{x}) = \alpha K(x).
        \end{equation}\end{linenomath*}
        By Equation \eqref{eq:ndK}, we see that  $\text{Supp} (\tilde{K})  = [-1, 1] $ and that 
        \begin{linenomath*}\begin{equation}
            \int_{-1}^{1} \tilde{K}(\tilde{x}) d\tilde{x} =  \int_{-1}^{1}  \alpha K(\alpha \tilde{x})  d\tilde{x} =\int_{-\alpha}^{\alpha} K(x) dx =1.
        \end{equation}\end{linenomath*}
    By \eqref{eq:nd} and \eqref{eq:ndK}, it follows that the convolution product becomes
	\begin{linenomath*}\begin{equation}\label{eq:ndconv}
	\begin{aligned}
	    K \ast u_i (x) &=  \int_{-\alpha}^{\alpha} K(x-y) u_i(y) dy \\
	    &= \int_{-1}^{1} \frac{1}{\alpha} \tilde{K}(\tilde{x}-\tilde{y}) \frac{1}{l}\tilde{u}_i(\tilde{y}) \alpha d\tilde{y} \\
	    &= \frac{1}{l} \tilde{K}\ast_{\sim} \tilde{u}_i (\tilde{x}),
	  \end{aligned}
	  \end{equation}\end{linenomath*}
	  where $\ast_{\sim}$ denotes the convolution operator in the rescaled spatial coordinate.
	  
	By substituting Equations \eqref{eq:nd}, \eqref{eq:ndK} and \eqref{eq:ndconv} in Equations \eqref{eq:system1}, we obtain the following non-dimensional system
	\begin{linenomath*}\begin{equation}\label{eq:system1ND}
		\begin{aligned}
			&\partial_ {\tilde{t}} \tilde{u}_1= \partial_{\tilde{x}\tilde{x}} \tilde{u}_1+ \frac{\gamma}{l D}\partial_{\tilde{x}} \left(\tilde{u}_1  \partial_{\tilde{x}}  (\tilde K \ast_{\sim} \tilde{u}_2) \right), \\
			&\partial_ {\tilde{t}} \tilde{u}_2= \partial_{\tilde{x}\tilde{x}} \tilde{u}_2+ \frac{\gamma}{l D}\partial_{\tilde{x}} \left(\tilde{u}_2  \partial_{\tilde{x}}  (\tilde K \ast_{\sim} \tilde{u}_1) \right),
		\end{aligned}
	\end{equation}\end{linenomath*}
    where $\tilde{x}\in \left[-\frac{l}{2\alpha}, \frac{l}{2\alpha}\right]$. By the relations in Equation \eqref{eq:nd}, the boundary conditions now read as:
	\begin{linenomath*}\begin{equation}\label{eq:bc_tilde}
		\tilde{u}_i\left(-\frac{l}{2\alpha},\tilde{t}\right) = \tilde{u}_i\left(\frac{l}{2\alpha},\tilde{t}\right), \,	\partial_{\tilde{x}}\tilde{u}_i\left(-\frac{l}{2\alpha},\tilde{t}\right) = \partial_{\tilde{x}} \tilde{u}_i\left(\frac{l}{2\alpha},\tilde{t}\right), \, \forall i \in \{1, \dots, N\} \text{ and } \tilde{t}\geq0.
	\end{equation}\end{linenomath*}
	The boundary conditions (Equation \eqref{eq:bc_tilde}) imply that the total mass of each population $\tilde{u}_i$ is conserved in time. Therefore, for $ i=1,2 $ and all $\tilde{t} \geq0$, the following identities hold
	\begin{linenomath*}\begin{equation}\label{eq:int_condnd}
		\int_{-l/2\alpha}^{l/2\alpha} \tilde{u}_i (0,\tilde{t}) \text{d} \tilde{x}= \int_{-l/2\alpha}^{l/2\alpha} \tilde{u}_i(\tilde{x},\tilde{t}) \text{d} \tilde{x}= \int_{-l/2}^{l/2} \frac{l}{\alpha} u_i(x,t) dx =  \frac{l}{\alpha}p_i,
	\end{equation}\end{linenomath*}
	where the second equality uses the identities in Equation \eqref{eq:nd} and the third equality uses Equation \eqref{eq:int_cond}.
        By Equation \eqref{eq:int_condnd} it follows that the non-dimensional system in \eqref{eq:system1ND} has a unique equilibrium point given by
        \begin{linenomath*}\begin{equation}\label{eq:equilibriumtilde}
             \mathbf{\bar{\tilde{u}}}:=(\bar{\tilde{u}}_1, \bar{\tilde{u}}_2)=\left(p_1,p_2\right).
        \end{equation}\end{linenomath*}
    To simplify the notation, we define $\tilde{\gamma}:=\frac{\gamma}{l D}$ and $L:=\frac{l}{\alpha}$, and by dropping the tildes, the non-dimensional system \eqref{eq:system1ND} reads as
	\begin{linenomath*}\begin{equation}\label{eq:ndsystem}
		\begin{aligned}
	        &\partial_t u_1= \partial_{xx} u_1+ \gamma\partial_{x} \left(u_1  \partial_{x}  (K \ast u_2) \right), \\
			&\partial_t u_2= \partial_{xx} u_2+ \gamma\partial_{x} \left(u_2 \partial_{x}  (K \ast u_1) \right),
		\end{aligned}
	\end{equation}\end{linenomath*}
	where $x\in \left[-\frac{L}{2}, \frac{L}{2}\right]$. The boundary conditions for System \eqref{eq:ndsystem} read as:
	\begin{linenomath*}\begin{equation}\label{eq:nd_bc}
		{u}_i\left(-\frac{L}{2},{t}\right) = {u}_i\left(\frac{L}{2},{t}\right), \,	\partial_{{x}}{u}_i\left(-\frac{L}{2},{t}\right) = \partial_{{x}} {u}_i\left(\frac{L}{2},{t}\right), \, \forall i \in \{1, \dots, N\} \text{ and } {t}\geq0.
	\end{equation}\end{linenomath*}

\subsection{{Linear stability analysis}}
    We now perform a linear stability analysis of system \eqref{eq:ndsystem} about the equilibrium point 
	\begin{linenomath*}\begin{equation}\label{eq:equilibrium}
		\mathbf{\bar{u}}=(\bar{u}_1, \bar{u}_2)=(p_1,p_2),
	\end{equation}\end{linenomath*}
    (see Equation \eqref{eq:equilibriumtilde}). To this end, we consider a perturbation of the homogeneous solution \eqref{eq:equilibrium} of the following form
	\begin{linenomath*}\begin{equation}\label{eq:w}
	  	\mathbf{w}=  \begin{pmatrix}
		u_1 - \bar{u}_1 \\u_2-\bar{u}_2
		\end{pmatrix}= \mathbf{u}^{(0)} e^{\lambda t + i q {x}},
	\end{equation}\end{linenomath*}	
    subject to boundary conditions \eqref{eq:nd_bc}, where $\mathbf{u}^{(0)} $ is a constant vector, $ \lambda \in \mathbb{R} $ is the growth rate and $ q $ is the wavenumber of the perturbation. By substituting Equation \eqref{eq:w} into Equation \eqref{eq:ndsystem} and neglecting nonlinear terms, we obtain the following eigenvalue problem
    \begin{linenomath*}\begin{linenomath*}\begin{equation}
		\lambda(q)\mathbf{w} = \mathcal{L}(q)\mathbf{w},
    \end{equation}\end{linenomath*}\end{linenomath*}
    where
    \begin{linenomath*}\begin{linenomath*}\begin{equation}\label{eq:L}
		\mathcal{L}(q)= -q^2 \begin{bmatrix}
			1  & \gamma \bar{u}_1 \hat{K}(q)  \\
			&  &  &  \\
			\gamma\bar{u}_2 \hat{K}(q) & 1
		\end{bmatrix},
    \end{equation}\end{linenomath*}\end{linenomath*}
    and
    \begin{linenomath*}\begin{equation}\label{eq:hatK}
        \begin{aligned}
	     \hat{K}(q):=& \int_{-1}^{1} K({x}) e^{-i  q {x}} \text{d} {x}  = \int_{-1}^{1} K({x}) \cos( q x)\text{d} {x},
        \end{aligned}
    \end{equation}\end{linenomath*}    
    where the second equality uses the fact that $K(x)$ is an even function and then $K(x)\sin(qx)$ is an odd function.

    The eigenvalues of the matrix $\mathcal{L}$ \eqref{eq:L} read 
    \begin{linenomath*}\begin{equation}\label{eq:lambdapm}
        \lambda^{\pm}(q) := - q^2 (1 \pm  \gamma \lvert\hat{K}(q)\rvert\sqrt{\bar{u}_1\bar{u}_2}),
    \end{equation}\end{linenomath*}
   and govern the evolution of the perturbation $\mathbf{w}$ (Equation \eqref{eq:w}). If $\gamma=0$ then $\lambda^{\pm}(q)\leq0$. By continuity, if $\gamma$ is arbitrarily small,  $\lambda^{\pm}(q)\leq0$ for all wavenumbers $q$, and the equilibrium point $\Bar{\mathbf{u}}$ (Equation \eqref{eq:equilibrium}) is linearly stable. As $|\gamma|$  increases, either $\lambda^{+}(q)$ or $\lambda^{-}(q)$ becomes positive for some values of $q$ and, consequently, the equilibrium $\bar{\mathbf{u}}$ becomes unstable.
   
   The wavenumbers $q$ must be chosen in such a way that the periodic boundary conditions in Equation \eqref{eq:nd_bc} are satisfied, and thus we have a discrete set of admissible wavenumbers given by
	\begin{linenomath*}\begin{equation}\label{eq:qm}
		I=\left\{q_m:= \frac{2 \pi m}{L}, \text{ with }  m \in \mathbb{Z}_{\geq 0}\right\}.
	\end{equation}\end{linenomath*}
   
   The equilibrium $\bar{\mathbf{u}}$ (Equation \eqref{eq:equilibrium}) is unstable when $\lambda^{\pm}(q_m)>0$ for some $m \in \mathbb{Z}_{\geq 0}$. Note that $\lambda^\pm(q_0)=0$ so the system never becomes unstable at wavenumber $q_0$. For $m>0$, if $\hat{K}(q_m)\neq 0$, we denote by $\gamma_m^{\pm}$ the instability thresholds of the wavenumber $q_m$, which are defined as
    \begin{linenomath*}\begin{equation}\label{eq:gammapm}
         \gamma_m^{\pm}= \pm \frac{1}{\lvert\hat{K}(q_m)\rvert\sqrt{\bar{u}_1\bar{u}_2}}, \, m \in \mathbb{Z}_{> 0}.
    \end{equation}\end{linenomath*} 
  Therefore the equilibrium $\bar{\mathbf{u}}$ (Equation \eqref{eq:equilibrium}) is unstable when
    \begin{linenomath*}\begin{equation}\label{eq:inst_cond}
           \gamma < \gamma_m^-  \quad  \text{ or } \quad \gamma> \gamma_m^+, \quad \text{ for some } m \in \mathbb{Z}_{>0}.
    \end{equation}\end{linenomath*}

    In the following section, we will perform a weakly nonlinear analysis to study the evolution of the perturbation $\mathbf{w}$ when the equilibrium $\bar{\mathbf{u}}$ becomes linearly unstable.
    We will adopt $\gamma$ as bifurcation parameter and denote by $q_c$ the first admissible wavenumber that is destabilized as $\lvert \gamma \rvert$ is increased. By Equation \eqref{eq:gammapm}, we note the critical wavenumber $q_c$ is defined as
    \begin{linenomath*}
        \begin{equation}
            q_c= \arg \max_{q_m \in I} \lvert \hat{K}(q_m) \rvert,
        \end{equation}
    \end{linenomath*}
    where the set $I$ is defined in \eqref{eq:qm}. We also underline that $q_c$ depends on the choice of kernel $K$ and may not be unique. We will denote by $\gamma_c^{\pm}$ the corresponding bifurcation thresholds, that is 
    \begin{linenomath*}\begin{equation}\label{eq:gamma_c}
    \gamma_c^+= \frac{1}{\lvert\hat{K}(q_c)\rvert\sqrt{\bar{u}_1\bar{u}_2}}\quad  \text{ and } \quad \gamma_c^- =- \frac{1}{\lvert\hat{K}(q_c)\rvert\sqrt{\bar{u}_1\bar{u}_2}}.
    \end{equation}\end{linenomath*}   

    \section{Amplitude equations}\label{sec:wnl_analysis}
    In this section we  perform a weakly nonlinear analysis based on the method of multiple scales. Close to the threshold of instability, that is in the weakly non-linear regime, we will use an expansion technique to recover an approximated solution, characterized by a slowly varying amplitude, and the equations governing the amplitude of the solution. Through the analysis of these equations (usually referred to as amplitude equations), we recover the amplitude and stability of the stationary solutions.
    
    The idea behind the multiple scale method comes from the observation that, just above an instability threshold, a nonlinear state is given by a superposition of modes whose wavenumbers $q$ lie in a narrow band $q^-\leq q \leq q^+$ (see \cite{CG09} Cap 6). The resulting nonlinear state is a solution governed by one or more unstable modes and characterized by an amplitude that varies slowly in space, due to the superposition of modes with almost identical wavenumbers. Also, the amplitude evolves slowly in time because, close to the onset of instability, all growth rates are small.

    Generally just beyond a bifurcation threshold, if the band of unstable wavenumbers $[q^-,q^+]$ around $q_c$ has width $O(\varepsilon)$, where $\varepsilon \ll 1$, the positive growth rates are $O(\varepsilon^2)$. Therefore, the solution evolves as
    \begin{linenomath*}\begin{equation}
        \begin{aligned}
             \mathbf{u}(x,t) & \sim \mathbf{\bar{u}} + \tilde{A}(X,T) e^{i q_c x} + \tilde{A^*}(X,T) e^{-i q_c x} ,
        \end{aligned}       
    \end{equation}\end{linenomath*}
    where $X=\varepsilon x$ is a long spatial scale, $T=\varepsilon^2 t$ is a slow temporal scale, $\tilde{A}(X,T)$ is a complex function and denotes the slow modulation of the critical mode $e^{i q_c x}$, and $\tilde{A^*}$ is the complex conjugate of $\tilde{A}$. Also, in the limit of $\varepsilon \rightarrow 0$, this solution must satisfy the boundary conditions in Equation \eqref{eq:nd_bc}. 
    
    However, in systems with a conservation law, so that $\lambda(0)=0$, long-scale modes evolve on long timescales, and must be included in the analysis (see also \cite{MC00}). Therefore solutions to System \eqref{eq:ndsystem}-\eqref{eq:nd_bc} evolve as

    \begin{linenomath*}\begin{equation}\label{eq:ansatz}
        \mathbf{u}(x,t)=\mathbf{\bar{u}} + \tilde{A}(X,T) e^{i q_c x} + \tilde{A^*}(X,T) e^{-i q_c x} + \tilde{B}(X,T) ,
    \end{equation}\end{linenomath*}
   where $\tilde{B}(X,T)$ is a real function and denotes the slow modulation of the mode corresponding to the zero wavenumber, $q=0$. 

    Recall that the homogeneous steady state is linearly stable for $\gamma_c^-<\gamma<\gamma_c^+$, and becomes unstable for $\gamma<\gamma_c^{-}$ or $\gamma>\gamma_c^+$.
    In the following Theorem, we derive an approximation of the solutions close to the instability thresholds ($\gamma \approx \gamma_c^{+}$ or $\gamma \approx \gamma_c^{-}$) and the equations governing the amplitude of the solutions. Since the analysis is broadly the same, we do not distinguish between $\gamma_c^{+}$ and $\gamma_c^{-}$ and use $\gamma_c$ to denote both the thresholds. 
    This Theorem also shows that the ansatz in Equation \eqref{eq:ansatz} correctly describes solutions in the weakly nonlinear regime.
    
    \begin{Theorem}\label{t:theorem}
    Let $ \varepsilon:= \sqrt{\lvert \frac{\gamma-\gamma_c}{\gamma_c} \rvert}$. When $\varepsilon \ll 1$, solutions to system \eqref{eq:ndsystem} have the following form
    \begin{linenomath*}\begin{equation}
    \begin{aligned}
        &u_1= \bar{u}_1 + \varepsilon \rho_1 ( A e^{i q_c x}  +  A^* e^{-i q_c x}) + \varepsilon^2[\psi_1 (A^2 e^{2i q_c x} + A^*{^2} e^{-2i q_c x}) + B] + O(\varepsilon^3),\\
        &u_2= \bar{u}_2 + \varepsilon \rho_2 ( A e^{i q_c x}  +  A^* e^{-i q_c x}) + \varepsilon^2[\psi_2 (A^2 e^{2i q_c x} + A^*{^2} e^{-2i q_c x}) + B] + O(\varepsilon^3).
    \end{aligned}
    \end{equation}\end{linenomath*}
    Here, $(\bar{u}_1,\bar{u}_2)$ is the homogeneous steady state \eqref{eq:equilibrium}, and $\rho_1$, $\rho_2$, $\psi_1$, $\psi_2$ are constants defined as
    \begin{linenomath*}\begin{equation}
    \begin{aligned}\label{rhos}
         &\rho_1=1, \qquad \rho_2=- \frac{1}{\gamma_c \bar{u}_1 \hat{K}(q_c)}, \\
         &\psi_1=  \frac{1}{2 \bar{u}_1} \frac{1-\gamma_c \bar{u}_1 \hat{K}(2q_c)}{1-\gamma_c^2 \bar{u}_1 \bar{u}_2 \hat{K}^2(2q_c)}, \qquad \psi_2 = \frac{1}{2 \bar{u}_1} \frac{1-\gamma_c \bar{u}_2 \hat{K}(2q_c)}{1-\gamma_c^2 \bar{u}_1 \bar{u}_2 \hat{K}^2(2q_c)}.
    \end{aligned}
    \end{equation}\end{linenomath*}
Also, $A(X,T)$ and $B(X,T)$ are governed by the following equations
    \begin{enumerate}
         \item If $\bar{u}_1 \neq \bar{u}_2$,
            \begin{linenomath*}\begin{equation}
		\begin{aligned}\label{eq:A2}
			&A_T= \sigma A - \Lambda \lvert A \rvert^2 A,\\
                &B=0,\end{aligned}.
            \end{equation}\end{linenomath*}    
        \item   If $\bar{u}_1=\bar{u}_2$,
            \begin{linenomath*}\begin{equation}\label{eq:AB2}
		\begin{aligned}
			&A_T= \sigma A - \Lambda \lvert A \rvert^2 A + \nu A B,\\
			&B_T= \mu B_{XX}-\eta( \lvert A \rvert^2)_{XX},
            \end{aligned}
	   \end{equation}\end{linenomath*}   
    \end{enumerate}
   where the coefficients $ \sigma $, $ \Lambda $, $ \nu $, $ \mu $ and $\eta$  are defined as
   \begin{linenomath*}\begin{equation}\label{eq:coefficients_ampl_eq}
       \begin{aligned}
           \sigma &= -q_c^2, \text{ if } \gamma_c^-<\gamma<\gamma_c^+ \text{ (stable regime)}, \quad  \sigma = q_c^2, \text{ if }  \gamma<\gamma_c^{-} \text{ or }\gamma>\gamma_c^+\text{ (unstable regime)},\\
           \Lambda &=\frac{1}{2} q_c^2 \gamma_c [2 \hat{K}(2 q_c)(\psi_1+\psi_2)-\hat{K}(q_c)(\psi_1 \rho_2 +\psi_2 a_2)], \\
           \nu &= \frac{q_c^2}{\bar{u}_1}, \qquad \mu=1+\gamma_c \bar{u}_1 \hat{K}(0), \qquad \eta=\frac{1}{\bar{u}_1}.
       \end{aligned}
   \end{equation}\end{linenomath*}
   Finally, $A^*$ denotes the complex conjugate of $A$.
    \end{Theorem}
    \begin{Proof} 
 Recall the definition of $\mathbf{w}$ from Equation \eqref{eq:w}.
	Separating the linear part from the non linear part, System \eqref{eq:ndsystem} can be rewritten as
	\begin{linenomath*}\begin{equation}\label{eq:weq}
		\partial_t \mathbf{w} =  \partial_{xx}\mathcal{L}^{\gamma} [\mathbf{w}] +  \partial_x \mathcal{Q}^{\gamma} [\mathbf{w},  \partial_x (K \ast \mathbf{w})],
	\end{equation}\end{linenomath*}
	where the actions of linear operator $ \mathcal{L}^{\gamma} $ and the non-linear operator $\mathcal{Q}^{\gamma}  $ on the vectors $ \mathbf{r} = (r_1, r_2)^T $ and $ \mathbf{s} = (s_1, s_2)^T $ are defined as
	\begin{linenomath*}\begin{equation}\label{eq:LQ}
		\mathcal{L}^{\gamma} \left[ \mathbf{r}\right]=\begin{pmatrix}
			1 && \gamma \bar{u}_1 K \ast  \\
			\gamma \bar{u}_2 K \ast && 1
		\end{pmatrix} \begin{pmatrix}
		    r_1 \\ r_2
		\end{pmatrix}, \qquad
		\mathcal{Q}^{\gamma}\left[ \mathbf{r}, \mathbf{s}\right] =
		\gamma 	\begin{pmatrix}
			r_1   s_2 \\ r_2 s_1
		\end{pmatrix}.
	\end{equation}\end{linenomath*}
        Choosing $\gamma$ such that $ \gamma - \gamma_c \sim \varepsilon^2 $, we write the following expansion
	\begin{linenomath*}\begin{equation}\label{eq:expgamma}
		\gamma=\gamma_c + \varepsilon^2 \gamma^{(2)}.
	\end{equation}\end{linenomath*}
    From the definition of $\varepsilon$, it follows that either $\gamma^{(2)}=  \gamma_c$ or $\gamma^{(2)}=  -\gamma_c$. In particular, $\gamma^{(2)}=-\gamma_c$ in the stable regime ($\gamma_c^-<\gamma<\gamma_c^+$), while $\gamma^{(2)}=\gamma_c$ in the unstable regime ($\gamma<\gamma_c^{-}$ or $\gamma>\gamma_c^+$).
    
    We then employ the method of multiple scales and adopt a long spatial scale $X=\varepsilon x$ and multiple temporal scales $T_1, T_2, \dots$ such that
    \begin{linenomath*}
        \begin{equation}
            t= \frac{T_1}{\varepsilon}+\frac{T_2}{\varepsilon^2}+ \cdots \,.
        \end{equation}
    \end{linenomath*}
    As $\varepsilon \rightarrow 0$, temporal and spatial derivatives decouple as
	\begin{linenomath*}\begin{equation}\label{eq:tx_der}
		\partial_t \rightarrow \partial_t +\varepsilon \partial_{T_1 }+ \varepsilon^2 \partial_{T_2}, \qquad \partial_x \rightarrow \partial_x + \varepsilon \partial_X.
	\end{equation}\end{linenomath*} 
    We employ a regular asymptotic expansion of $\mathbf{w}$ in terms of  $\varepsilon$
	\begin{linenomath*}\begin{equation}\label{eq:expw}
			\mathbf{w}= \varepsilon \mathbf{w}_1 +\varepsilon ^2 \mathbf{w}_2 +\varepsilon^3 \mathbf{w}_3 + \cdots,
	\end{equation}\end{linenomath*}
    where
        \begin{linenomath*}\begin{equation}\label{eq:wj}
            \mathbf{w}_j = \sum_{m=-\infty}^\infty \mathbf{w}_{j m}(X,T_1,T_2) e^{i q_m x}, \, \text{ for } j=1,2, \dots
        \end{equation}\end{linenomath*}
    and must satisfy the boundary conditions in Equations \eqref{eq:nd_bc}.
    
	By Equations \eqref{eq:expgamma} and \eqref{eq:expw}, we see that the operators $ \mathcal{L}^{\gamma} $ and $\mathcal{Q}^{\gamma}  $ in \eqref{eq:LQ} decouple in orders of $\varepsilon$ as
	\begin{linenomath*}\begin{equation}\label{eq:LQexp}
	\begin{aligned}
		\mathcal{L}^{\gamma} \left[ \mathbf{r}\right]&=\begin{pmatrix}
			1 && (\gamma_c+\varepsilon^2 \gamma^{(2)}) \bar{u}_1 K \ast  \\
			(\gamma_c+\varepsilon^2 \gamma^{(2)}) \bar{u}_2 K \ast && 1
		\end{pmatrix} \begin{pmatrix}
		    r_1 \\ r_2
		\end{pmatrix}\\
		&= \mathcal{L}^{\gamma_c} \left[ \mathbf{r}\right]+\varepsilon^2\begin{pmatrix}
		0	 &&  \gamma^{(2)} \bar{u}_1 K \ast  \\
		\gamma^{(2)} \bar{u}_2 K \ast && 0
		\end{pmatrix} \begin{pmatrix}
		    r_1 \\ r_2
		\end{pmatrix},\\ \\
		\mathcal{Q}^{\gamma}\left[ \mathbf{r}, \mathbf{s}\right] &=
		(\gamma_c+\varepsilon^2 \gamma^{(2)}) 	\begin{pmatrix}
			r_1   s_2 \\ r_2 s_1
		\end{pmatrix}= \mathcal{Q}^{\gamma_c}\left[ \mathbf{r}, \mathbf{s}\right] + \varepsilon^2 \mathcal{Q}^{\gamma^{(2)}}\left[ \mathbf{r}, \mathbf{s}\right].
		\end{aligned}
	\end{equation}\end{linenomath*}
	By substituting Equations \eqref{eq:expw}, \eqref{eq:expgamma}, \eqref{eq:tx_der} and \eqref{eq:LQexp} into Equation \eqref{eq:weq}, we obtain
	\begin{linenomath*}\begin{equation}
	    \begin{aligned}
	        &\varepsilon^2 \partial_{T_1} \mathbf{w}_1+\varepsilon^3 \partial_{T_2} \mathbf{w}_1 +\varepsilon^3 \partial_{T_1} \mathbf{w}_2 +\varepsilon^4 \partial_{T_2} \mathbf{w}_2=\\
	        & \quad (\partial_{xx}+2 \varepsilon \partial_{xX}+\varepsilon^2 \partial_{XX}) \mathcal{L}^{\gamma_c}[\varepsilon \mathbf{w}_1 +\varepsilon ^2 \mathbf{w}_2 +\varepsilon^3 \mathbf{w}_3+\varepsilon^4 \mathbf{w}_4] \\
	        & \quad+ \varepsilon^2 (\partial_{xx}+2 \varepsilon \partial_{xX}+\varepsilon^2 \partial_{XX}) \begin{pmatrix}
		0	 &&  \gamma^{(2)} \bar{u}_1 K \ast  \\
		\gamma^{(2)} \bar{u}_2 K \ast && 0
		\end{pmatrix} (\varepsilon \mathbf{w}_1 +\varepsilon ^2 \mathbf{w}_2 +\varepsilon^3 \mathbf{w}_3) \\
	        & \quad +(\partial_{x}+\varepsilon \partial_{X}) \mathcal{Q}^{\gamma_c}[(\varepsilon \mathbf{w}_1 +\varepsilon ^2 \mathbf{w}_2 +\varepsilon^3 \mathbf{w}_3),(\partial_{x}+\varepsilon \partial_{X})(K \ast (\varepsilon \mathbf{w}_1 +\varepsilon ^2 \mathbf{w}_2 +\varepsilon^3 \mathbf{w}_3))] \\
	        & \quad + \varepsilon^2(\partial_{x}+\varepsilon \partial_{X}) \mathcal{Q}^{\gamma^{}(2)}[(\varepsilon \mathbf{w}_1 +\varepsilon ^2 \mathbf{w}_2 +\varepsilon^3 \mathbf{w}_3),(\partial_{x}+\varepsilon \partial_{X})(K \ast (\varepsilon \mathbf{w}_1 +\varepsilon ^2 \mathbf{w}_2 +\varepsilon^3 \mathbf{w}_3))] + O(\varepsilon^5).
	    \end{aligned}
	\end{equation}\end{linenomath*}
    Next we collect the terms at each order of $ \varepsilon $ and obtain a sequence of equations for each $ \mathbf{w}_i $.	At order $ \varepsilon $, we obtain the homogeneous linear problem $ \partial_{xx} \mathcal{L}^{\gamma_c}[\mathbf{w}_1] =0$, where the function $\mathbf{w}_1$, has
    the form as in \eqref{eq:wj}. Therefore, we have:
    \begin{linenomath*}
        \begin{equation}\label{eq:sol1}
        \begin{aligned}
             \partial_{xx} \mathcal{L}^{\gamma_c}[\mathbf{w}_1] &= \partial_{xx}\sum_{m=-\infty}^{\infty} 
             \begin{pmatrix}
                 1 & \gamma_c \bar{u}_1 K \ast   \\ \gamma_c \bar{u}_2 K \ast  & 1
              \end{pmatrix}  \mathbf{w}_{1m} e^{i q_m x}
              \\&= \partial_{xx} \begin{pmatrix}
                 1 & \gamma_c \bar{u}_1 \hat{K}(q_m)   \\ \gamma_c \bar{u}_2 \hat{K}(q_m) & 1
                 \end{pmatrix} \mathbf{w}_{1m} e^{i q_m x}
             \\&=-\sum_{m=-\infty}^{\infty}q_m^2
             \begin{pmatrix}
                 1 & \gamma_c \bar{u}_1 \hat{K}(q_m)   \\ \gamma_c \bar{u}_2 \hat{K}(q_m) & 1
                 \end{pmatrix} \mathbf{w}_{1m} e^{i q_m x}\\&=0
                 \end{aligned}
        \end{equation}
    \end{linenomath*}
    where the second equality uses 
    \begin{linenomath*}
        \begin{equation}\label{eq:hat_K}
            K \ast e^{i q_m x} = \int_{-1}^{1} K(y) e^{i q_m (x-y)}dy = \int_{-1}^{1} K(y) e^{-i q_m y}dy e^{i q_m x} = \hat{K}(q_m) e^{i q_m x}.
        \end{equation}
    \end{linenomath*}
    with $ \hat{K} $ defined in \eqref{eq:hatK}.
    The fourth equality in Equation \eqref{eq:sol1} is satisfied if and only if 
    \begin{linenomath*}
        \begin{equation}\label{eq:solm}
        q_m^2    \begin{pmatrix}
                 1 & \gamma_c \bar{u}_1 \hat{K}(q_m)   \\ \gamma_c \bar{u}_2 \hat{K}(q_m) & 1
                 \end{pmatrix} \mathbf{w}_{1m} e^{i q_m x}=0, \, \text{for all } m \in \mathbb{Z}.
        \end{equation}
    \end{linenomath*}
    Non-trivial solutions to Equation \eqref{eq:solm} exist when either the determinant of the matrix is zero or $q_m=0$. Recalling the definition of $q_m $ \eqref{eq:qm} and $\gamma_c$ \eqref{eq:gamma_c}, we see that non-trivial solutions exist only for $q_m=q_0$ and $q_m=q_c$.

    Therefore, the function $\mathbf{w}_1$ that satisfies this linear problem $ \partial_{xx} \mathcal{L}^{\gamma_c}[\mathbf{w}_1] =0$ is
	\begin{linenomath*}\begin{equation}\label{eq:w1}
		\mathbf{w}_1 = \bm{\rho}_{0} A_0(X,T_1,T_2)+\bm{\rho}\left(A(X,T_1,T_2) e^{i q_c x} +A^*(X,T_1,T_2) e^{-i q_c x} \right) ,
	\end{equation}\end{linenomath*}
	where $A_0(X,T_1,T_2)$ is a real function, $A(X,T_1,T_2)$ is a complex function,  $A^* $ denotes the complex conjugate of $ A $, and
	\begin{linenomath*}\begin{equation}\label{eq:rhos}
		\bm{\rho}_0 = \begin{pmatrix} \rho_{01} \\ \rho_{02}	\end{pmatrix}, \qquad \bm{\rho} =  \begin{pmatrix} \rho_1 \\ \rho_2	\end{pmatrix}
	\end{equation}\end{linenomath*}
	are constant vectors. First, notice that $  \partial_{xx} \mathcal{L}^{\gamma_c}[\bm{\rho}_{0} A_0(X,T_1,T_2)] =0 $, for any $ \bm{\rho}_{0} $ and $ A_0(X,T_1,T_2) $. Also, in order to satisfy $ \partial_{xx} \mathcal{L}^{\gamma_c}[\mathbf{w}_1] =0$, the vector $ \bm{\rho} $ must be such that
	\begin{linenomath*}\begin{equation}\label{eq:rhokern}
		\bm{\rho} \in \text{Ker} \begin{pmatrix}
			1 && \gamma_c \bar{u}_1 \hat{K}(q_c)\\
		\gamma_c \bar{u}_2 \hat{K}(q_c) && 1
		\end{pmatrix},
	\end{equation}\end{linenomath*}
	where $ \hat{K} $ is defined in \eqref{eq:hatK}.
	Since $ \gamma_c \hat{K}(q_c)  \sqrt{\bar{u}_1\bar{u}_2}=\pm1 $ (see Equation \eqref{eq:gamma_c}),  $	\bm{\rho} $ can be defined up to a constant. We shall choose the following normalization
	\begin{linenomath*}\begin{equation}\label{eq:rho}
		\bm{\rho} = \begin{pmatrix} 1 \\ \rho_2 \end{pmatrix}, \text{ where } \rho_2 :=- \frac{1}{\gamma_c \bar{u}_1 \hat{K}(q_c) }.
	\end{equation}\end{linenomath*}
	At this stage, the amplitudes $ A(X,T_1,T_2) $ and $ A_0(X,T_1,T_2) $, and the vector $ \bm{\rho}_0 $ are still unknown.
	
	At order $ \varepsilon^2 $ we obtain the following problem
	\begin{linenomath*}\begin{equation}\label{eq:eps2}
			\partial_{xx} \mathcal{L}^{\gamma_c} [\mathbf{w}_2]  = \mathbf{F},
	\end{equation}\end{linenomath*}
	with
	\begin{linenomath*}\begin{equation}\label{eq:F}
		\begin{aligned}		
		\mathbf{F}=	& -2 \partial_{xX}\mathcal{L}^{\gamma_c} [\mathbf{w}_1] -\partial_x \mathcal{Q}^{\gamma_c} [\mathbf{w}_1, \partial_x (K \ast \mathbf{w}_1)] +\partial_{T_1} \mathbf{w}_1\\
			=&- 2iq_c \begin{pmatrix}
			1 && \gamma_c \bar{u}_1 \hat{K}(q_c)\\
		\gamma_c \bar{u}_2 \hat{K}(q_c) && 1
		\end{pmatrix} \begin{pmatrix}
		    1 \\ \rho_2
		\end{pmatrix}	(A_X e^{i q_c x} - A^*_X e^{-i q_c x} ) \\
			&+2 q_c^2 \rho_2 \gamma_c \hat{K}(q_c) \begin{pmatrix} 1 \\ 1	\end{pmatrix}	(A^2 e^{2i q_c x} + A^*{^2} e^{-2i q_c x} ) \\
			&+ q_c^2  \gamma_c \hat{K}(q_c)  \begin{pmatrix} \rho_{01} \rho_2 \\ \rho_{02}	\end{pmatrix} A_0 (A e^{i q_c x} + A^* e^{-i q_c x} ) \\ & +\bm{\rho}_{0} \partial_{T_1} A_{0}+\bm{\rho}\left(\partial_{T_1} A e^{i q_c x} +\partial_{T_1} A^* e^{-i q_c x} \right) 
		\\	=& \, 2 q_c^2 \rho_2 \gamma_c \hat{K}(q_c) \begin{pmatrix} 1 \\ 1	\end{pmatrix}	(A^2 e^{2i q_c x} + A^*{^2} e^{-2i q_c x} ) + q_c^2  \gamma_c \hat{K}(q_c)  \begin{pmatrix} \rho_{01} \rho_2 \\ \rho_{02}	\end{pmatrix} A_0 (A e^{i q_c x} +A^* e^{-i q_c x} )\\
			=& \, -\frac{2}{\bar{u}_1} q_c^2  \begin{pmatrix} 1 \\ 1	\end{pmatrix}	(A^2 e^{2i q_c x} + A^*{^2} e^{-2i q_c x} ) + q_c^2  \gamma_c \hat{K}(q_c)  \begin{pmatrix} \rho_{01} \rho_2 \\ \rho_{02}	\end{pmatrix} A_0 (A e^{i q_c x} +A^* e^{-i q_c x} )
   \\ & +\bm{\rho}_{0} \partial_{T_1} A_{0}+\bm{\rho}\left(\partial_{T_1} A e^{i q_c x} +\partial_{T_1} A^* e^{-i q_c x} \right) ,
		\end{aligned} 
	\end{equation}\end{linenomath*}
	where the second equality uses Equation \eqref{eq:hat_K}, the third equality is true because, by Equation \eqref{eq:rhokern}, the term on the second line is equal to zero, and the fourth equality uses the definition of $\rho$ (Equation \eqref{eq:rho}).
	
	By the Fredholm Alternative Theorem, Equation \eqref{eq:eps2} admits a solution if and only if for any $\mathbf{a} \in L^2(-L/2, L/2)$ such that
	\begin{linenomath*}\begin{equation}\label{eq:a}
	    \mathbf{a} \in \text{Ker} \{(\partial_{xx}\mathcal{L}^{\gamma_c})^{T}\} = \text{Ker} \left\{\partial_{xx}\begin{pmatrix}
			1 && \gamma_c \bar{u}_2 K \ast  \\
			\gamma_c \bar{u}_1 K \ast && 1
		\end{pmatrix}\right\},
	\end{equation}\end{linenomath*}
	the equality $\langle \mathbf{F}, \mathbf{a} \rangle=0$ is satisfied, where $\langle \cdot, \cdot \rangle$ denotes the scalar product in $L^2(-L/2, L/2)$.
	Notice that any  $\mathbf{a} \neq \mathbf{0}$ satisfying the condition in \eqref{eq:a} is a constant multiple of
	\begin{linenomath*}\begin{equation}\label{eq:aest}
	    \mathbf{a} = \begin{pmatrix}
	        1 \\ a_2
	    \end{pmatrix} (e^{i q_c x}+e^{-i q_c x}), \text{ with } a_2 :=- \frac{1}{\gamma_c \bar{u}_2 \hat{K}(q_c) }.
	\end{equation}\end{linenomath*}
	Therefore Equation \eqref{eq:eps2} only has a solution when $ \rho_{01}=\rho_{02} = 0$ and $\partial_{T_1} A=0$, that is the amplitude $A$ does not depend on $T_1$. From now on, we will denote $T_2$ by $T$ for simplicity and write $A(X,T)$ instead of $A(X,T_2)$. 
 
    Therefore, the linear problem in Equation \eqref{eq:eps2} reduces to
	\begin{linenomath*}\begin{equation}\label{eq:eps2_2}
		\partial_{xx} \mathcal{L}^{\gamma_c} [\mathbf{w}_2] =-\frac{2}{\bar{u}_1} q_c^2 \begin{pmatrix} 1 \\ 1	\end{pmatrix}	(A^2(X,T) e^{2i q_c x} + A^*{^2}(X,T) e^{-2i q_c x} ).
	\end{equation}\end{linenomath*}	
	Finally, by Equation \eqref{eq:eps2_2} it follows that the function $ \mathbf{w}_2 $, having the form as in \eqref{eq:wj}, is given by
	\begin{linenomath*}\begin{equation}\label{eq:w2}
		\mathbf{w}_2 =\bm{\psi}_0 B_0(X,T)+ \bm{\psi} (A^2(X,T) e^{2i q_c x} + A^*{^2}(X,T) e^{-2i q_c x}),
	\end{equation}\end{linenomath*}
	where $B_0(X,T)$ is a real function and
	\begin{linenomath*}\begin{equation}\label{eq:psis}
		\bm{\psi}_0 = \begin{pmatrix} \psi_{01} \\ \psi_{02}	\end{pmatrix}, \qquad \bm{\psi} = \begin{pmatrix} \psi_1 \\ \psi_2	\end{pmatrix}
	\end{equation}\end{linenomath*}
	are constant vectors. Notice that $  \partial_{xx} \mathcal{L}^{\gamma_c}[\bm{\psi}_0 B_0(X,T)] =0 $, for any $ \bm{\psi}_{0} $ and $ B_0(X,T) $. Substituting Equation \eqref{eq:w2} into Equation \eqref{eq:eps2_2} and solving for $\bm{\psi} $ we obtain
		\begin{linenomath*}\begin{equation}\label{eq:psi}
		\begin{aligned}
		   &\psi_1 = \frac{1}{2 \bar{u}_1} \frac{1-\gamma_c \bar{u}_1 \hat{K}(2q_c)}{1-\gamma_c^2 \bar{u}_1 \bar{u}_2 \hat{K}^2(2q_c)}, \\
		   &\psi_2 = \frac{1}{2 \bar{u}_1} \frac{1-\gamma_c \bar{u}_2 \hat{K}(2q_c)}{1-\gamma_c^2 \bar{u}_1 \bar{u}_2 \hat{K}^2(2q_c)},
		\end{aligned}
	\end{equation}\end{linenomath*}
        whilst $ A(X,T) $, $ B_0(X,T) $, and $ \bm{\psi}_0 $ remain unknown.
        
	At order $ \epsilon^3 $, we find the following problem
	\begin{linenomath*}\begin{equation}\label{eq:eps3}
	    \partial_{xx} \mathcal{L}^{\gamma_c} [\mathbf{w}_3] = \mathbf{G},
	\end{equation}\end{linenomath*}
	where
	\begin{linenomath*}\begin{equation}\label{eq:G}
	\begin{aligned}
	    \mathbf{G} = & \, \partial_T \mathbf{w}_1 - 2\partial_{xX}  \mathcal{L}^{\gamma_c} [\mathbf{w}_2] - \partial_{XX}  \mathcal{L}^{\gamma_c} [\mathbf{w}_1] - \begin{pmatrix}0 && \gamma^{(2)} \bar{u}_1 K \ast \\  \gamma^{(2)} \bar{u}_2 K \ast && 0 \end{pmatrix} \partial_{xx}  \mathbf{w}_1 \\
			& - \partial_{x} \mathcal{Q}^{\gamma_c}[\mathbf{w}_1, \partial_{x}(K \ast \mathbf{w}_2)] -  \partial_{x} \mathcal{Q}^{\gamma_c}[\mathbf{w}_2, \partial_{x}(K \ast \mathbf{w}_1)]  \\
			& - \partial_{x} \mathcal{Q}^{\gamma_c}[\mathbf{w}_1, \partial_{X}(K \ast \mathbf{w}_1)]  -  \partial_{X} \mathcal{Q}^{\gamma_c}[\mathbf{w}_1, \partial_{x}(K \ast \mathbf{w}_1)] \\
			= & \, (A_T e^{i q_c x} + A_T^* e^{-i q_c x}) \bm{\rho}+ 8 i q_c \begin{pmatrix}
			1 && \gamma_c \bar{u}_1 \hat{K}(2q_c)\\
		\gamma_c \bar{u}_2 \hat{K}(2q_c) && 1
		\end{pmatrix} \bm{\psi} (A A_X e^{2iq_c x}-A^{*} A^{*}_X e^{-2iq_c x})\\
		&- \begin{pmatrix}
			1 && \gamma_c \bar{u}_1 \hat{K}(q_c)\\
		\gamma_c \bar{u}_2 \hat{K}(q_c) && 1
		\end{pmatrix} \bm{\rho}	(A_{XX} e^{i q_c x} + A^*_{XX} e^{-i q_c x} )
			\\ & + q_c^2 \gamma^{(2)} \hat{K}(q_c) \begin{pmatrix} \rho_2 \bar{u}_1 \\ \bar{u}_2	\end{pmatrix}(A e^{i q_c x} + A^* e^{-i q_c x} ) 
			\\&+  q_c^2\gamma_c   \left( 2  \hat{K}(2 q_c) \begin{pmatrix} \psi_2 \\ \psi_1 \rho_2	\end{pmatrix}-\hat{K}(q_c)  \begin{pmatrix} \psi_1 \rho_2 \\ \psi_2	\end{pmatrix}\right)\lvert A \rvert^2 (A e^{i q_c x} + A^* e^{-i q_c x} ) \\
			&  +q_c^2 \gamma_c \hat{K}(q_c) \begin{pmatrix} \psi_{01}\rho_2 \\ \psi_{02}	\end{pmatrix} (A e^{i q_c x} + A^* e^{-i q_c x} )  B_0\\
			&  +3  q_c^2\gamma_c \left( 2  \hat{K}(2 q_c) \begin{pmatrix}  \psi_2  \\  \psi_1 \rho_2	\end{pmatrix}+\hat{K}(q_c)  \begin{pmatrix} \psi_1 \rho_2 \\ \psi_2	\end{pmatrix}\right) (A^3 e^{3i q_c x} + A^*{^3} e^{-3i q_c x} ) \\
			& - 4 i q_c \hat{K}(q_c) \rho_2 \begin{pmatrix}  1 \\ 1	\end{pmatrix} (A A_X e^{2iq_c x}-A^{*} A^{*}_X e^{-2iq_c x}).
	 \end{aligned}
	\end{equation}\end{linenomath*}
	 By Equation \eqref{eq:rhokern}, it follows that the third term of the second equality of Equation \eqref{eq:eps3} is the null vector. In order to simplify the notation, we rewrite Equation \eqref{eq:G} as:
	 \begin{linenomath*}\begin{equation}\label{eq:GG}
        \begin{aligned}
	     \mathbf{G} = & (A_T \bm{\rho} + A \mathbf{G}_1 + A B_0 \mathbf{G}_1^{(2)} + |A|^2 A \mathbf{G}_1^{(3)})e^{i q_c x} + \mathbf{G}_2 (A^2)_{X} e^{2i q_c x} + \mathbf{G}_3 A^3 e^{3i q_c x} 
      \\ & + (A^*_T \bm{\rho} + A^* \mathbf{G}_1 + A^* B_0 \mathbf{G}_1^{(2)} + |A^*|^2 A^* \mathbf{G}_1^{(3)})e^{-i q_c x} + \mathbf{G}_2 (A^{*^2})_{X} e^{-2i q_c x} + \mathbf{G}_3 A^{*^3} e^{-3i q_c x}.
      \end{aligned}
	 \end{equation}\end{linenomath*}

    The linear problem in Equation \eqref{eq:eps3} admits a solution if and only the Fredholm condition $\langle \mathbf{G}, \mathbf{a} \rangle =0 $ is satisfied, where $\mathbf{a}$ is defined in Equation \eqref{eq:a}. Note that the terms $ \mathbf{G}_2 (A^2)_{X} e^{2i q_c x}+\mathbf{G}_2 (A^{*^2})_{X} e^{-2i q_c x}$ and $\mathbf{G}_3 A^3 e^{3i q_c x} +\mathbf{G}_3 A^{*^3} e^{-3i q_c x}$ are hortogonal to $\mathbf{a}$. 
    Therefore, the Fredholm condition $\langle \mathbf{G}, \mathbf{a} \rangle =0 $ for Equation \eqref{eq:eps3} gives the following amplitude equation
	\begin{linenomath*}\begin{equation}\label{eq:A}
		A_T= \sigma A - \Lambda \lvert A \rvert^2 A + \delta A B_0,
	\end{equation}\end{linenomath*}
	where
	\begin{linenomath*}\begin{equation}\label{eq:Acoeff}
	\begin{aligned}
	    \sigma &= - \frac{\langle \mathbf{G}_1, \mathbf{a}\rangle}{\langle \bm{\rho}, \mathbf{a}\rangle} = q_c^2 \frac{\gamma^{(2)}}{\gamma_c}
	    \\ \Lambda &= \frac{\langle \mathbf{G}_1^{(3)}, \mathbf{a}\rangle}{\langle \bm{\rho}, \mathbf{a}\rangle} = \frac{1}{2} q_c^2 \gamma_c (2 \hat{K}(2 q_c)(\psi_1+\psi_2)-\hat{K}(q_c)(\psi_1 \rho_2 +\psi_2 a_2)) \\
	    \delta &= - \frac{\langle \mathbf{G}_1^{(2)}, \mathbf{a}\rangle}{\langle \bm{\rho}, \mathbf{a}\rangle} = \frac{1}{2} q_c^2 \left(\frac{\psi_{01}}{\bar{u}_1}+\frac{\psi_{02}}{\bar{u}_2}\right).
	 \end{aligned}
	\end{equation}\end{linenomath*}
   
    At order $ \epsilon^4 $, we have the following problem
    \begin{linenomath*}\begin{equation}\label{eq:eps4}
		\begin{aligned}
		\partial_{xx} \mathcal{L}^{\gamma_c} [\mathbf{w}_4] = & \, 	\partial_T \mathbf{w}_2  - 2\partial_{xX}  \mathcal{L}^{\gamma_c} [\mathbf{w}_3] - \partial_{XX}  \mathcal{L}^{\gamma_c} [\mathbf{w}_2] \\
			&-\gamma^{(2)} \hat{K}(2q_c) \begin{pmatrix}0 && \bar{u}_1 \\  \bar{u}_2 && 0 \end{pmatrix} \partial_{xx}  \mathbf{w}_2 - 2\gamma^{(2)} \hat{K}(q_c)\begin{pmatrix}0 && \bar{u}_1 \\ \bar{u}_2 && 0 \end{pmatrix} \partial_{xX}  \mathbf{w}_1 \\
			& - \partial_{x} \mathcal{Q}^{\gamma_c}[\mathbf{w}_1, \partial_{x}(K \ast \mathbf{w}_3)] -  \partial_{x} \mathcal{Q}^{\gamma_c}[\mathbf{w}_3, \partial_{x}(K \ast \mathbf{w}_1)] -  \partial_{x} \mathcal{Q}^{\gamma_c}[\mathbf{w}_2, \partial_{x}(K \ast \mathbf{w}_2)] \\
			& - \partial_{x} \mathcal{Q}^{\gamma_c}[\mathbf{w}_1, \partial_{X}(K \ast \mathbf{w}_2)]  - \partial_{x} \mathcal{Q}^{\gamma_c}[\mathbf{w}_2, \partial_{X}(K \ast \mathbf{w}_1)] \\
			& -  \partial_{X} \mathcal{Q}^{\gamma_c}[\mathbf{w}_1, \partial_{x}(K \ast \mathbf{w}_2)] -  \partial_{X} \mathcal{Q}^{\gamma_c}[\mathbf{w}_2, \partial_{x}(K \ast \mathbf{w}_1)] \\
			& - \partial_{X} \mathcal{Q}^{\gamma_c}[\mathbf{w}_1, \partial_{X}(K \ast \mathbf{w}_1)]  -  \partial_{x} \mathcal{Q}^{\gamma^{(2)}}[\mathbf{w}_1, \partial_{x}(K \ast \mathbf{w}_1)] \\
			= & \, \begin{pmatrix} \psi_{01} \\ \psi_{02} \end{pmatrix} (B_0)_T - \begin{pmatrix} \psi_{01}+\gamma_c \bar{u}_1 \hat{K}(0) \psi_{02}\\\gamma_c \bar{u}_2 \hat{K}(0) \psi_{01} +\psi_{02}	\end{pmatrix} (B_0)_{XX} + \frac{1}{\bar{u}_1} \begin{pmatrix} 1 \\ 1 \end{pmatrix} (|A|^2)_{XX} \\
			& + \sum_{h=1}^{3} \mathbf{r}_h e^{i h q_c x} + \mbox{c.c.}
		\end{aligned}
	\end{equation}\end{linenomath*}
    Since the function $\mathbf{w}_4$ is as in \eqref{eq:wj}, in Equation \eqref{eq:eps4} all terms independent of $x$ must be equal to zero, that is
    \begin{linenomath*}\begin{equation}\label{eq:dtB0}
        \begin{pmatrix} \psi_{01} \\ \psi_{02} \end{pmatrix} (B_0)_T=  \begin{pmatrix} \psi_{01}+\gamma_c \bar{u}_1 \hat{K}(0) \psi_{02}\\\gamma_c \bar{u}_2 \hat{K}(0) \psi_{01} +\psi_{02}	\end{pmatrix} (B_0)_{XX} - \frac{1}{\bar{u}_1} \begin{pmatrix} 1 \\ 1 \end{pmatrix} (|A|^2)_{XX}.
    \end{equation}\end{linenomath*}
    When $\bar{u}_1=\bar{u}_2$, we can choose $\psi_{01}=\psi_{02}$ and, by setting $B:=\psi_{01} B_0$, we obtain the following amplitude equations
	\begin{linenomath*}\begin{equation}\label{eq:AB}
            \begin{aligned}
                A_T&= \sigma A - \Lambda \lvert A \rvert^2 A + \nu A B, \\
		      B_T&= \mu B_{XX} -  \eta (\lvert A \rvert^2)_{XX},  
            \end{aligned}
	\end{equation}\end{linenomath*}
	where 
	\begin{linenomath*}\begin{equation}\label{eq:Bcoeff}
		\nu= \frac{q_c^2}{\bar{u}_1}, \qquad \mu=1+\gamma_c \bar{u}_1 \hat{K}(0), \qquad \eta=\frac{1}{\bar{u}_1},
	\end{equation}\end{linenomath*}
    and $\sigma$ and $\Lambda$ are given in Equation \eqref{eq:Acoeff}. Notice that $\nu = \delta/\psi_{01} $ (see Equation \eqref{eq:Acoeff}), with $\psi_{01}=\psi_{02}$ and $\bar{u}_1=\bar{u}_2$. On the other hand, if $\bar{u}_1 \neq \bar{u}_2$, Equation \eqref{eq:dtB0} is satisfied when $\psi_{01}=\psi_{02}=0$ and $(|A|^2)_{XX}=0$.
    \end{Proof}	
\newline
\subsection{{Small amplitude solutions}} 
   The stationary solutions of the amplitude equations in \eqref{eq:A2} and \eqref{eq:AB2} correspond to steady states of system \eqref{eq:ndsystem}. 
    Notice that if $ B=0 $, Equation \eqref{eq:AB2} reduces to Equation \eqref{eq:A2}, which is a Stuart-Landau equation. If $ \Lambda>0$, system \eqref{eq:ndsystem} undergoes a
    supercritical bifurcation, while if $ \Lambda<0$ the system undergoes a subcritical bifurcation (\cite{TB90}). 

    In the supercritical regime, as the homogeneous steady state becomes unstable, stationary small amplitude patterns emerge and correspond to solutions of Equation \eqref{eq:A2} with $ A=a_0 e^{i \phi} $, where $ \phi \in \mathbb{R} $ is the phase of the pattern and the amplitude $ a_0 $ is real and must satisfy $ a_0^2 = \sigma/\Lambda $. These small amplitude solutions are always stable (\cite{TB90}). 
    
    Analogously, stationary small amplitude patterns correspond to solutions of Equation \eqref{eq:AB2} with $ A=a_0 e^{i \phi} $ and $ B=0 $, where $ \phi \in \mathbb{R} $ and $ a_0^2 = \sigma/\Lambda $. However, in this case the stationary patterns might be destabilized by large-scale modes (\cite{CH93}). In the following Proposition we will derive a stability condition for these stationary solutions.
    \begin{Proposition}\label{pr:stability}
        Suppose $\bar{u}_1=\bar{u}_2$. If $ \sigma>0$ and $ \Lambda>0$ then small amplitude patterns to System \eqref{eq:ndsystem} exist. These solutions are unstable if the following condition holds
        \begin{linenomath*}\begin{equation}\label{eq:Gamma}
            	\Gamma:= \frac{\Lambda \mu}{\eta \nu}-1<0,
        \end{equation}\end{linenomath*}
        where $\sigma$, $\Lambda$, $\mu$, $\eta$ and $\nu$ are given in Theorem \ref{t:theorem}.
        \end{Proposition}
    \begin{Proof}
     By Theorem \ref{t:theorem}, if $\bar{u}_1=\bar{u}_2$, the amplitude of the stationary solutions to System \eqref{eq:ndsystem} is governed by Equation \eqref{eq:AB2}.  When $ \sigma>0$ and $ \Lambda>0$, stationary small amplitude patterns exist and correspond to solutions of \eqref{eq:AB2} with $ A=a_0 e^{i \phi} $ and $ B=0 $, where $ \phi \in \mathbb{R} $ and $ a_0^2 = \sigma/\Lambda $.
    To study the stability of this stationary solution, we consider the following perturbation
	\begin{linenomath*}\begin{equation}\label{eq:ABpert}
		A(X,T)=(a_0 + a(X,T))e^{i \phi} , \qquad B(X,T)=b(X,T).
	\end{equation}\end{linenomath*}
    We substitute the perturbation \eqref{eq:ABpert} in Equations \eqref{eq:AB2}, and by linearizing in $ a $ and $ b $ we obtain:
	\begin{linenomath*}\begin{equation}\label{eq:ab_lin}
	\begin{aligned}
			a_T=& -\sigma (a  + a^*) +\nu a_0 b,\\
			b_T=& \mu b_{XX} - \eta a_0 (a_{XX} +a^*_{XX}).
	\end{aligned}
	\end{equation}\end{linenomath*}
	We consider a perturbation of the form
	\begin{linenomath*}\begin{equation}\label{eq:pert}
	    a(X,T)= e^{\bar{\lambda} T} (V e^{iQX} + W^*e^{-iQX}) \text{ and } b(X,T)= e^{\bar{\lambda}T} (U e^{iQX}  +U^*e^{-iQX}),
	\end{equation}\end{linenomath*}
    where $\bar{\lambda}$ is the growth rate of the perturbation, $U,V,W \in \mathbb{C}$ {and $Q\geq 0$ denotes a spatial mode}. Notice that $a$ is a complex perturbation, while $b$ is real. Upon substituting Equations \eqref{eq:pert} in Equations \eqref{eq:ab_lin}, we obtain the following eigenvalue problem
    \begin{linenomath*}\begin{equation}\label{eq:VWU}
		\bar{\lambda} \begin{pmatrix} V \\ W \\ U \end{pmatrix}= 
             \begin{pmatrix}
		    -\sigma & -\sigma & \nu a_0\\
                -\sigma & -\sigma & \nu a_0\\
                 \eta a_0 Q^2&  \eta a_0 Q^2 &  -\mu Q^2
	    \end{pmatrix} \begin{pmatrix}  V \\ W \\ U \end{pmatrix},      
    \end{equation}\end{linenomath*}
    from which we recover the growth rates
    \begin{linenomath*}\begin{equation}
        \bar{\lambda}_0(Q)=0, \qquad \bar{\lambda}^{\pm}(Q)=\frac{1}{2} \left(-\mu Q^2-2 \sigma \pm \sqrt{\mu ^2 Q^4+Q^2 \left(8 a_0^2 \eta  \nu -4 \mu  \sigma \right)+4 \sigma ^2} \right).
    \end{equation}\end{linenomath*}
    Recalling that $a_0^2=\sigma/\Lambda$, a simple calculation shows that $\bar{\lambda}^+(Q)>0$ if $ Q \neq 0$ and $\Gamma= \frac{\Lambda \mu}{\eta \nu}-1<0$.
    \end{Proof}
\newline
\newline

The analysis so far is valid for any non-negative, symmetric kernel $K$ satisfying Equation (\ref{eq:kint}). In the following section, we adopt the top-hat distribution and use the results obtained so far to recover the instability thresholds and to predict the shape of the emerging patterns.

For readers convenience, we conclude this section with Table \eqref{tab:table}, in which we have collected the main parameters involved in the study of stability and bifurcations, and included a brief description of their significance and properties.

\begin{table}[h]
    \centering
    \begin{tabular}{ |c|c|c|}
 \hline
 \multicolumn{3}{|c|}{\textbf{Parameter List and Description}} \\
 \hline
 \textbf{Parameter}  & \textbf{Description}    & \textbf{Properties}\\
 \hline
 $\gamma:= \gamma_{12}=\gamma_{21}$   & Inter-species interaction parameter & $\gamma>0:$ Mutual avoidance \\
   &  &  $\gamma<0:$ Mutual attraction \\
 \hline
 $L:= l/\alpha $ &   Length of the rescaled domain   &  $L>2$: ratio between the length \\
 & & of the domain $l$ and the sensing range $\alpha$\\
 \hline
 $\sigma$ (Eq. \eqref{eq:Acoeff}) & Linear Stuart-Landau coefficient & $\sigma<0$:  $\mathbf{\bar{u}}$ (Eq. \eqref{eq:equilibrium}) stable\\
 &  & $\sigma>0$: $\mathbf{\bar{u}}$ (Eq. \eqref{eq:equilibrium})  unstable \\
 \hline
 $\Lambda $ (Eq. \eqref{eq:Acoeff}) & Nonlinear Stuart-Landau coefficient & $\Lambda<0$: subcritical bifurcation \\
 &  & $\Lambda>0$: supercritical bifurcation\\
 \hline
  $\Gamma$ (Eq.\eqref{eq:Gamma})   & Stability coefficient  & $\Gamma<0$: unstable supercritical bifurcation \\
 & computed for $\Lambda>0$ & $\Gamma>0$: stable supercritical bifurcation\\
 \hline
\end{tabular}
    \caption{List and description of main parameters involved in the study of stability and bifurcations.}
    \label{tab:table}
\end{table}

\section{The top hat distribution}\label{sec:th}

In this section we analyze System \eqref{eq:system1} with
\begin{linenomath*}\begin{equation}\label{eq:top-hat}
	K(x)=K_{\alpha}(x):= 
	\begin{cases}
		\frac{1}{2 \alpha}, \quad  x\in [-\alpha, \alpha]\\
		0, \quad \text{ otherwise}
	\end{cases}.
\end{equation}\end{linenomath*}
The parameter $\alpha$, modelling the sensing radius of an organism, is such that $\alpha < l/2$, where $l$ is the length of the domain.  As in Section \ref{sec:lin_analysis}, we will work in dimensionless co-ordinates, so that our study system is given by Equations (\ref{eq:ndsystem}) and the dimensionless averaging kernel is
	\begin{linenomath*}\begin{equation}\label{eq:ndtop-hat}
	   K_1(x)=\begin{cases}
		\frac{1}{2}, \quad  x\in [-1, 1],\\
	0, \quad \text{ otherwise}.
	\end{cases}
	\end{equation}\end{linenomath*}
    

\subsection{Linear stability analysis}\label{sec:th_linear_analysis}
    Linear stability analysis of System \eqref{eq:ndsystem} around the equilibrium point $\mathbf{\bar{u}}=(p_1,p_2)$ (Equation \eqref{eq:equilibrium}), gives the following eigenvalues (see Equation \eqref{eq:lambdapm})
    \begin{linenomath*}\begin{equation}\label{eq:lambdapm_th}
        \lambda^{\pm}(q) := - q^2 (1 \pm  \gamma \lvert\hat{K}_1(q)\rvert\sqrt{\bar{u}_1\bar{u}_2}),
    \end{equation}\end{linenomath*}
    where
        \begin{linenomath*}\begin{equation}\label{eq:hat-kernel}
	     \hat{K}_{1}(q)= \int_{-1}^{1} K_{1}({x}) e^{-i  q {x}} \text{d} {x} = \begin{cases}
			\frac{\sin(q )}{q }, & \text{if } q\neq0 \\
		1, & \text{if } q=0
		\end{cases}.
        \end{equation}\end{linenomath*}
    Recall that the admissible wavenumbers are $q_m=2 \pi m/L$, with $m \in \mathbb{N}$.  
   
\begin{figure}[H]
        \centering
		\subfloat[Mutual avoidance ($ \gamma>0 $)]
		{\includegraphics[width=1\textwidth]{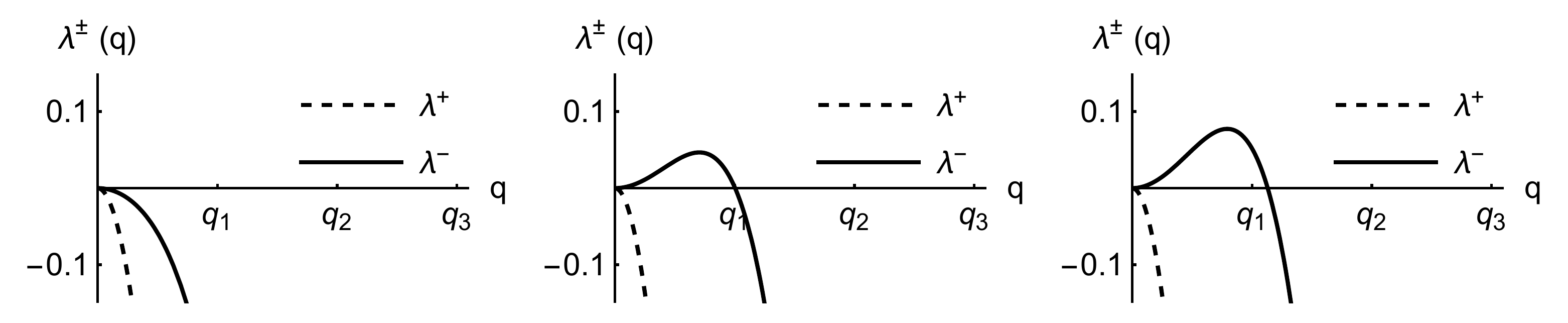}}\\
		\subfloat[Mutual attraction ($ \gamma<0 $)]
		{\includegraphics[width=1\textwidth]{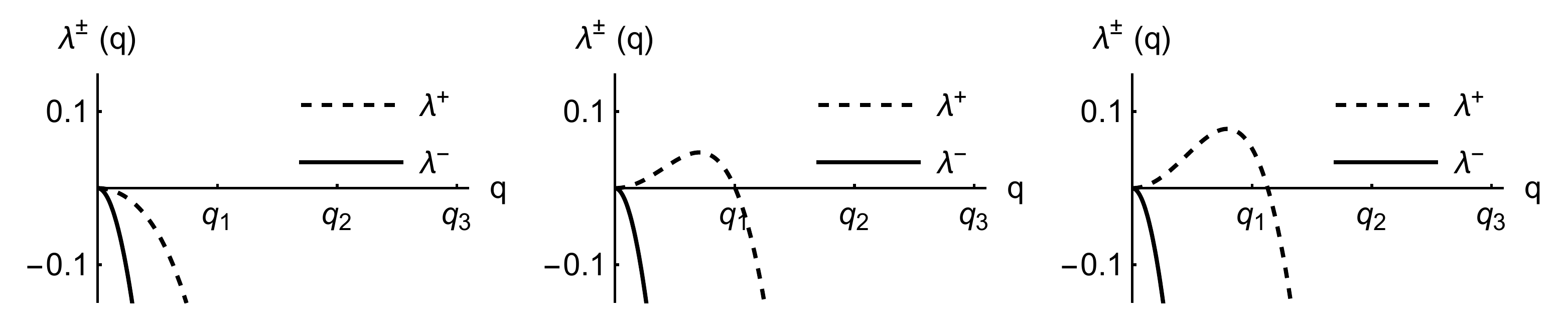}}
		\caption{Graphs of  the growth rates $ \lambda^{\pm}(q) $ (Equation \eqref{eq:lambdapm_th}), in the mutual avoidance (Panel (a)) and in the mutual attraction (Panel (b)) regime. 
            Panel (a) shows the graphs of $ \lambda^{\pm}(q) $ for increasing values of $\gamma>0$: $0<\gamma<\gamma_1^+$ (left); $\gamma=\gamma_1^+$ (center);  $\gamma>\gamma_1^+$ (right). Panel (b) shows the graphs of $ \lambda^{\pm}(q) $ for decreasing values of $\gamma<0$: $\gamma_1^-<\gamma<0$ (left); $\gamma=\gamma_1^-$ (center);  $\gamma<\gamma_1^-$ (right). As the magnitude of $\gamma$ increases, the first wavenumber destabilized is $q_1$ (Equation \eqref{eq:qc})}
		\label{fig:PlotEigen}
	\end{figure}

   Figure \ref{fig:PlotEigen} shows the graphs of $\lambda^\pm(q)$ (Equation \eqref{eq:lambdapm_th}) for different values of $\gamma$. Observe that the first wavenumber that is destabilized as $ \gamma $ is varied is
	\begin{linenomath*}\begin{equation}\label{eq:qc}
	   q_c=q_1 = \frac{2 \pi}{L}.
	\end{equation}\end{linenomath*}
    Since $L>2$, we have $\hat{K}_{1}(q_c)>0$, so the corresponding bifurcation thresholds, obtained by solving $\lambda^\pm(q_c)=0$, are 
	\begin{linenomath*}\begin{equation}\label{eq:gammac}
		\gamma_c^{\pm} = \gamma_1^{\pm} :=\pm \frac{1}{\hat{K}_{1}(q_c)\sqrt{\bar{u}_1\bar{u}_2}}. 
	\end{equation}\end{linenomath*}
    Since the equilibrium $\bar{\mathbf{u}}$  becomes unstable as $\lambda^{\pm}(q_c)>0$, the system undergoes an instability when
	\begin{linenomath*}\begin{equation}\label{eq:inst_cond_th}
		\gamma< \gamma_c^-=-\frac{1}{\hat{K}_{1}(q_c) \sqrt{\bar{u}_1\bar{u}_2}} \quad \text{ or } \quad \gamma > \gamma_c^+=\frac{1}{\hat{K}_{1}(q_c) \sqrt{\bar{u}_1\bar{u}_2}}.
	\end{equation}\end{linenomath*}
\subsection{Analysis of the amplitude equations and bifurcations}\label{sec:th_bifurcations}
    By Theorem \ref{t:theorem}, when $ \varepsilon= \sqrt{\lvert \frac{\gamma-\gamma_c}{\gamma_c} \rvert} \ll 1$ (where $\gamma_c=\gamma_c^{\pm}$), the solutions to System \eqref{eq:ndsystem} have the following form
    \begin{linenomath*}\begin{equation}\label{eq:u1u2}
    \begin{aligned}
        &u_1= \bar{u}_1 + \varepsilon \rho_1 ( A e^{i q_c x}  +  A^* e^{-i q_c x}) + \varepsilon^2(\psi_1 (A^2 e^{2i q_c x} + A^*{^2} e^{-2i q_c x}) + B) + O(\varepsilon^3),\\
        &u_2= \bar{u}_2 + \varepsilon \rho_2 ( A e^{i q_c x}  +  A^* e^{-i q_c x}) + \varepsilon^2(\psi_2 (A^2 e^{2i q_c x} + A^*{^2} e^{-2i q_c x}) + B) + O(\varepsilon^3).
    \end{aligned}
    \end{equation}\end{linenomath*}
    Recall from (\ref{rhos}) that the constants $\rho_1$, $\rho_2$ are defined as
    \begin{linenomath*}\begin{equation}
    \begin{aligned}
         &\rho_1=1, \qquad \rho_2=- \frac{1}{\gamma_c \bar{u}_1 \hat{K}_1(q_c)}.
    \end{aligned}
    \end{equation}\end{linenomath*}
    Note that in the mutual avoidance case ($ \gamma>0 $), $\gamma_c=\gamma_c^+>0$ and then $\rho_2<0$, which implies that $u_1$ and $u_2$ {show a spatial  oscillation that is  out of phase. On the other hand, in the mutual attraction regime ($ \gamma<0 $), $\gamma_c=\gamma_c^-<0$ and then $\rho_2>0$, which means that the spatial pattern for $u_1$ and $u_2$  are in phase.}

    Theorem \ref{t:theorem} also says that $A(X,T)$ and $B(X,T)$ are governed by the following equations
    \begin{enumerate}
         \item If $\bar{u}_1 \neq \bar{u}_2$,
            \begin{linenomath*}\begin{equation}
		\begin{aligned}\label{eq:A3}
			&A_T= \sigma A - \Lambda \lvert A \rvert^2 A,\\
                &B=0,\end{aligned}.
            \end{equation}\end{linenomath*}    
        \item   If $\bar{u}_1=\bar{u}_2$,
            \begin{linenomath*}\begin{equation}\label{eq:AB3}
		\begin{aligned}
			&A_T= \sigma A - \Lambda \lvert A \rvert^2 A + \nu A B,\\
			&B_T= \mu B_{XX}-\eta( \lvert A \rvert^2)_{XX},
            \end{aligned}
	   \end{equation}\end{linenomath*}   
    \end{enumerate}
   where the coefficients $ \sigma $, $ \Lambda $, $ \nu $, $ \mu $ and $\eta$  are defined in Equation \eqref{eq:coefficients_ampl_eq}
   
\begin{figure}[H]
        \centering
		\subfloat[Mutual avoidance ($ \gamma>0 $), $\bar u_1=\bar u_2$]
		{\includegraphics[width=1\textwidth]{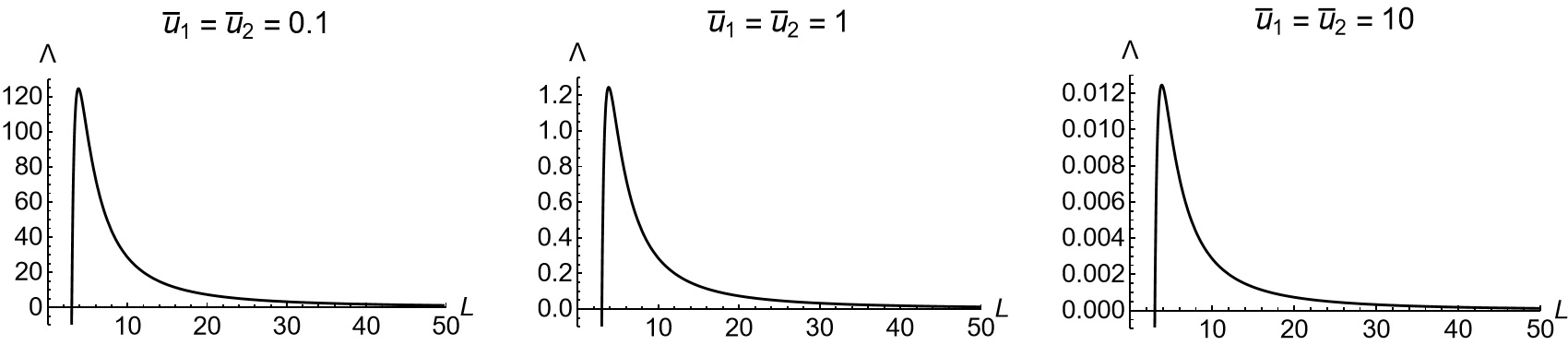}}\\
		\subfloat[Mutual avoidance ($ \gamma>0 $), $ \bar u_1 \neq \bar u_2$]
		{\includegraphics[width=1\textwidth]{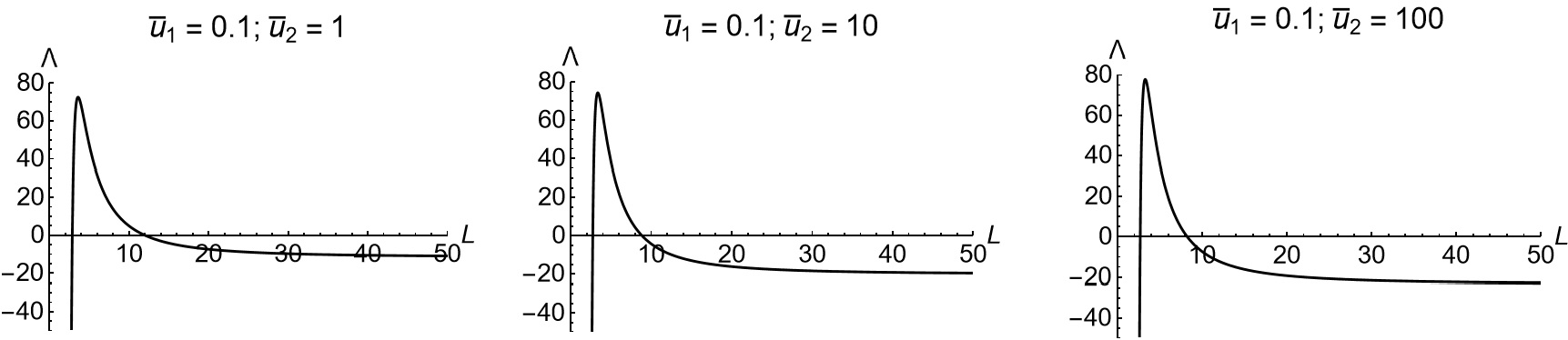}}\\
            \subfloat[Mutual attraction ($ \gamma<0 $), $\bar u_1=\bar u_2$]
		{\includegraphics[width=1\textwidth]{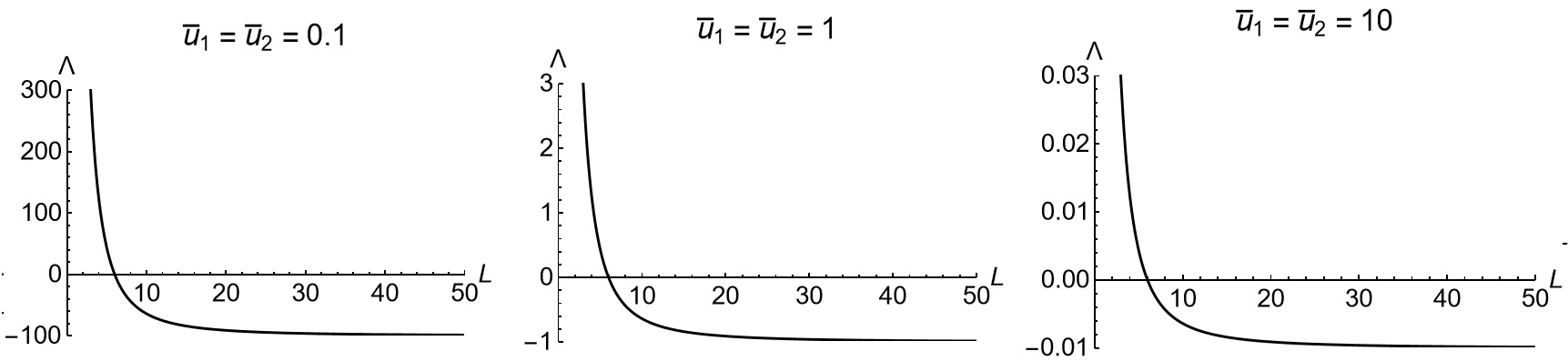}}\\
		\subfloat[Mutual attraction ($ \gamma<0 $), $\bar u_1 \neq \bar u_2$]
		{\includegraphics[width=1\textwidth]{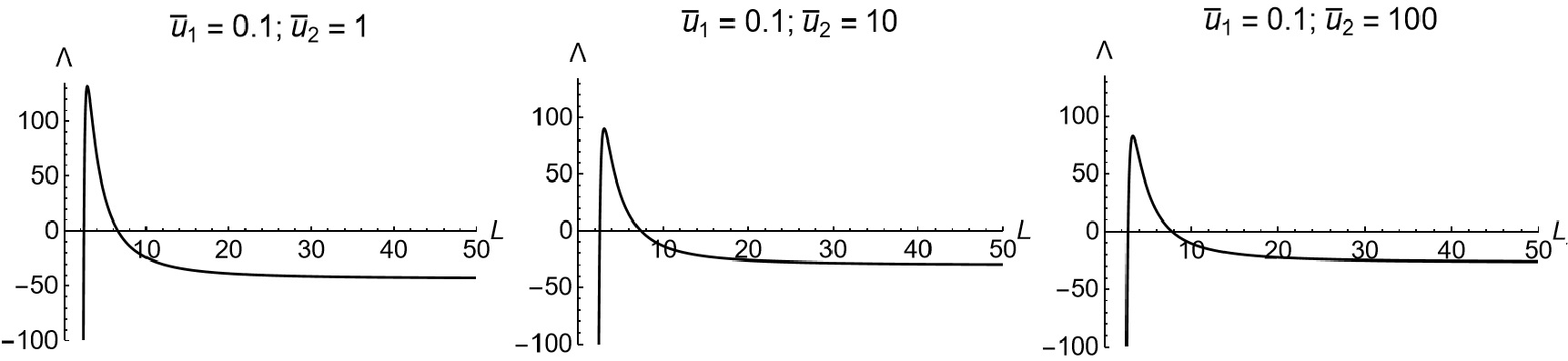}}
		\caption{Graphs of the nonlinear Stuart-Landau coefficient $\Lambda$ (Equation \eqref{eq:Acoeff}) versus the domain length $L$, in the mutual avoidance ($ \gamma>0 $) and in the mutual attraction ($ \gamma>0 $) regime, with $\bar u_1=\bar u_2$ and $ \bar u_1 \neq \bar u_2$. Positive values of $\Lambda$ correspond to supercritical bifurcations, negative values of $\Lambda$ correspond to subcritical bifurcations}
		\label{fig:L}
    \end{figure}

    As discussed in  Section \ref{sec:wnl_analysis}, the sign of $\Lambda$ determines the type of bifurcation: for $\Lambda>0$ the system exhibits a supercritical bifurcation, while for $\Lambda<0$ the system undergoes a subcritical bifurcation (see also Table \eqref{tab:table}). The sign of $\Lambda$ depends on $\bar{u}_1$, $\bar{u}_2$ and on the length of the domain, $L$ (see the definition of $\Lambda$ in Equation \eqref{eq:coefficients_ampl_eq}). Figure \ref{fig:L} shows the graphs of $\Lambda$ versus $L$, in the mutual avoidance ($\gamma>0$) and in the mutual attraction ($\gamma<0$) regime with $K=K_1$, for both $\bar u_1=\bar u_2$ and $ \bar u_1 \neq \bar u_2$.

    For $\gamma>0$, if $\bar{u}_1=\bar{u}_2$ then the qualitative behaviour of $\Lambda(L)$ remains unchanged as $\bar{u}_1=\bar{u}_2$ are varied.  In fact, Figure \ref{fig:L}(a) shows that for different values of $\bar{u}_1=\bar{u}_2$, $\Lambda(L)$ is negative (subcritical bifurcation) for $2<L<3$, while it is positive (supercritical bifurcation) for $L>3$. On the other hand, if $\bar{u}_1\neq\bar{u}_2$ (Figure \ref{fig:L}(b)), $\Lambda(L)$ is negative for $2<L<3$,  becomes positive for $L>3$, and then $\Lambda(L)$ becomes negative again for sufficiently large values of $L$ depending on the ratio $\bar{u}_1/\Bar{u}_2$.

    For $\gamma<0$, if $\bar{u}_1=\bar{u}_2$ , $\Lambda(L)$ is positive for $2<L<6$ and it becomes negative as $L>6$ (see Figure \ref{fig:L}(c)). The qualitative behaviour of  $\Lambda(L)$ does not change as $\bar{u}_1=\bar{u}_2$ are varied. However, if $\bar{u}_1\neq\bar{u}_2$ (Figure \eqref{fig:L} (d)) we observe the emergence of a subcritical regime for sufficiently small values of $L$ depending on the ratio $\bar{u}_1/\Bar{u}_2$.

    As shown in Section \ref{sec:wnl_analysis}, if $\Lambda(L)$ is positive then small amplitude patterns emerge from the homogeneous steady state beyond the bifurcation threshold. These solutions are always stable when $\bar{u}_1\neq\bar{u}_2$ but can be unstable when $\bar{u}_1=\bar{u}_2$. 
    
\begin{figure}[H]
        \centering
		\subfloat[Mutual avoidance ($ \gamma>0 $), $\bar u_1=\bar u_2$]
		{\includegraphics[width=0.42\textwidth]{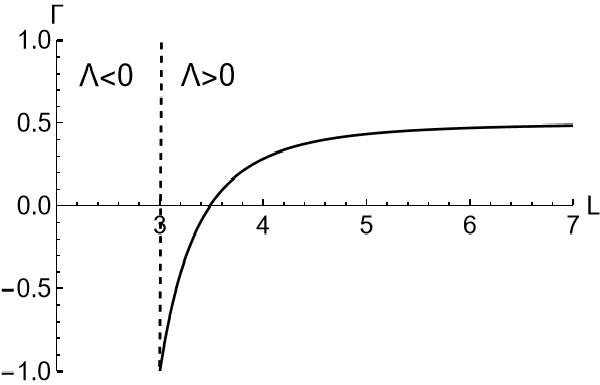}}
            \hspace{0.5cm}
            \subfloat[Mutual attraction ($ \gamma<0 $), $\bar u_1=\bar u_2$]
		{\includegraphics[width=0.42\textwidth]{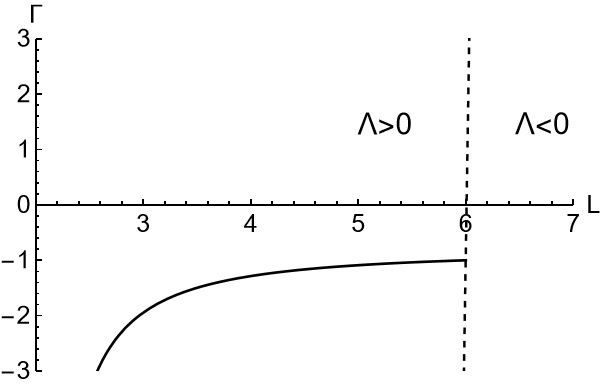}}
		\caption{Graphs of the stability coefficient  $\Gamma$ (Equation \eqref{eq:Gamma}) versus the domain length $L$, in the mutual avoidance (a), and in the mutual attraction regime (b). If $\Lambda>0$ and $\Gamma<0$, small amplitude patterns exist and are unstable, and if $\Lambda>0$ and $\Gamma>0$, small amplitude patterns exist and are stable}
		\label{fig:Gamma}
	\end{figure}

\begin{figure}[H]
        \centering
		\subfloat[Mutual avoidance ($ \gamma>0 $), $\bar u_1=\bar u_2$]
		{\includegraphics[width=0.42\textwidth]{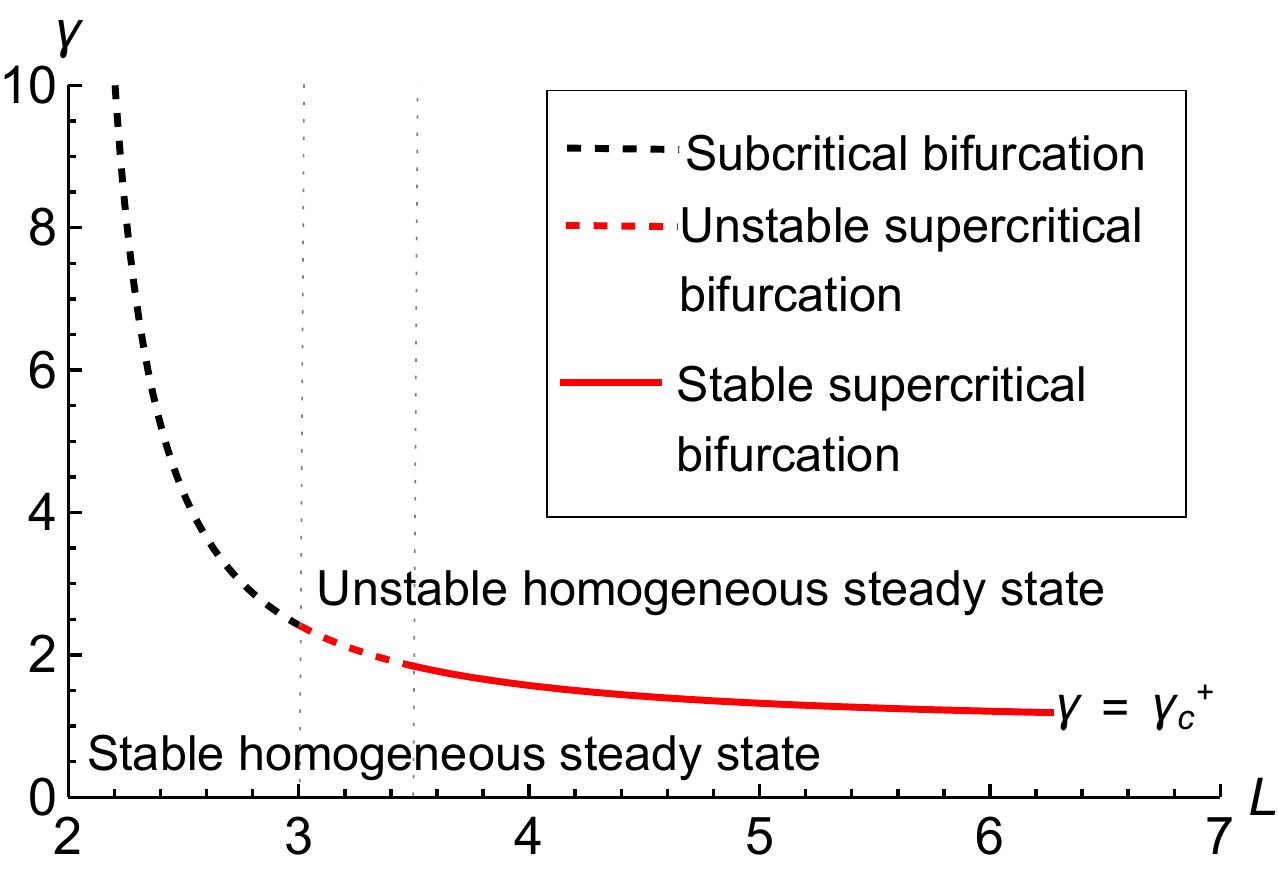}}\hspace{1cm}
            \subfloat[Mutual attraction ($ \gamma<0 $), $\bar u_1=\bar u_2$]
		{\includegraphics[width=0.42\textwidth]{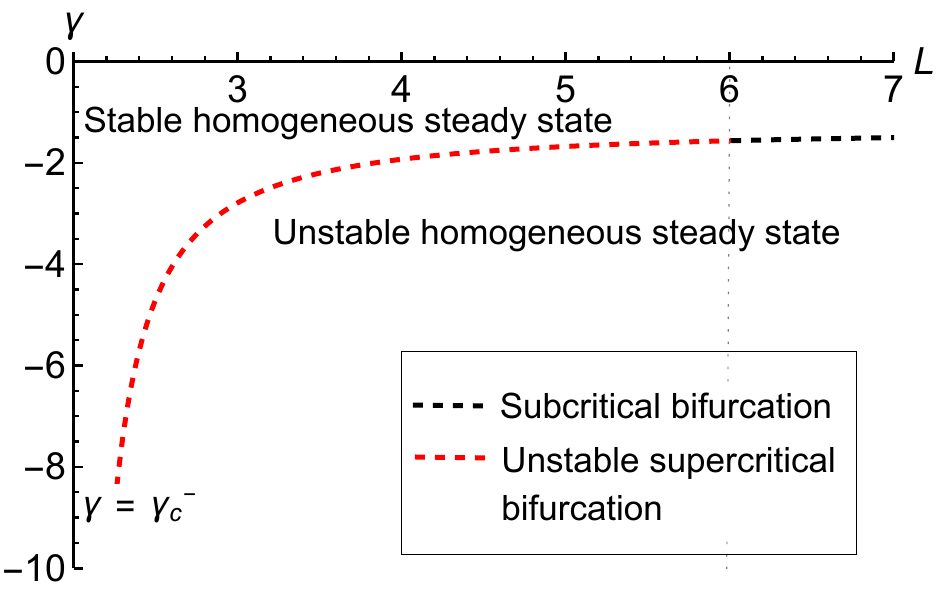}}
   		\caption{Graphs of the curves for the critical values of the density-dependent advection strength $\gamma=\gamma_c^+$ in (a) and $\gamma=\gamma_c^-$ in (b) (Equation \eqref{eq:gamma_c}) versus the domain length $L$. When the magnitude of $\gamma$ is small, the homogeneous steady state is linearly stable. As the magnitude of $\gamma$ increases, the system undergoes a bifurcation and the homogeneous steady state becomes unstable as $\gamma$ crosses $\gamma_c^{\pm}$. For $ \gamma>0 $ (a), when $L$ is small the system undergoes a subcritical bifurcation. As $L$ increases, the bifurcation becomes supercritical, and the emerging patterns will be unstable. As $L$ increases further, the system undergoes a supercritical bifurcation leading to the emergence of stable patterns. For $ \gamma<0 $ (b), when $L$ is small the system undergoes a supercritical bifurcation generating unstable small amplitude patterns. As $L$ increases, the bifurcation becomes subcritical}
        \label{fig:gammac_vs_L}
	\end{figure}

    Proposition \ref{pr:stability} shows that when $\bar{u}_1=\bar{u}_2$ the stability of small amplitude patterns is determined by the coefficients of the amplitude equations in \eqref{eq:AB3} and that, in particular, these solutions are unstable if $\Gamma=\frac{\Lambda \mu}{\eta \nu}-1<0$. By using the definitions of $ \Lambda $, $ \nu $, $ \mu $ and $\eta$ in Equation \eqref{eq:coefficients_ampl_eq}, we recover
    \begin{linenomath*}\begin{equation}\label{eq:Gamma_th}
       \Gamma=\frac{(1+\hat{K}_1(q_1))(2\hat{K}_1(2q_1)+\hat{K}_1(q_1) )}{2\hat{K}_1(q_1)(\hat{K}_1(2q_1)+\hat{K}_1(q_1))} -1.
    \end{equation}\end{linenomath*}
    Note that $\Gamma$ does not depend on $\bar{u}_1$. Indeed, since $q_1=2\pi/L$, it follows that $\Gamma$ depends only on $L$. In Figure \ref{fig:Gamma} we show the graphs of $\Gamma$ versus $L$ for $\gamma>0$ in (a), and $\gamma<0$ in (b). We also recall that we are analyzing the sign of $\Gamma$ in supercritical regimes ($\Lambda>0$), for this reason we plot the curve $\Gamma(L)$ only in those intervals in which $\Lambda>0$. The graph in Figure \ref{fig:Gamma}(a) shows that in the mutual avoidance case ($\gamma>0$), small amplitude patterns exist and are unstable for $3<L<3.5 $, and that they become stable as $L>3.5$. Figure \ref{fig:Gamma}(b) shows that in the mutual attraction scenario ($\gamma<0$), $\Gamma(L)$ is always negative and therefore small amplitude patterns are always unstable. These results are summarized in Figure \ref{fig:gammac_vs_L}.


    In summary, our analysis shows that the nature of the transition and the stability of the bifurcation patterns depend mainly on $L$. These results can be read and re-interpreted in terms of the parameters of the original system \eqref{eq:system1}, recalling that $L=\alpha/l$, where $\alpha$ is the sensing radius and $l$ is the length of the dimensional spatial domain. Therefore, the qualitative behaviour of the system under study strongly depends on the measure of the sensing radius compared on the length of the domain.

 \subsection{Numerical Simulations}\label{sec:th_simulations}
In this Section, we perform a numerical investigation of system \eqref{eq:ndsystem}. To solve numerically System \eqref{eq:ndsystem}, we use the spectral method and numerical schemes presented in \cite{GHLP22}. By employing a continuation technique, we recover numerical bifurcation diagrams which are compared with the bifurcation diagrams obtained via the weakly nonlinear analysis. We show that our weakly nonlinear analysis provides accurate approximations of stable steady-state solutions in supercritical stable regimes, as long as we stay close to the bifurcation threshold. We also analyse those bifurcations that generate unstable small amplitude patterns. In these cases, we numerically detect the existence of stable large amplitude solutions, which are not predicted by the weakly nonlinear analysis, {but which were predicted by an energy method in \cite{GHLP22}}. 

First, we analyze the scenarios depicted in Figures \ref{fig:L}(b) ($\gamma>0$) and (d) ($\gamma<0$), in which $\bar{u}_1\neq\bar{u}_2$. These figures show subcritical bifurcations for sufficiently small values of $L$, then a shift to a supercritical regime, as $L$ increases, and again a subcritical regime, as $L$ increases further. Recall that if $\bar{u}_1\neq\bar{u}_2$ then supercritical bifurcations always give rise to stable small amplitude solutions.

\begin{figure}[H]
		\centering
    $\gamma>0$, $\bar{u}_1\neq \bar{u}_2$
          \subfloat[$L=2.7$]
        {\includegraphics[width=0.33\textwidth]{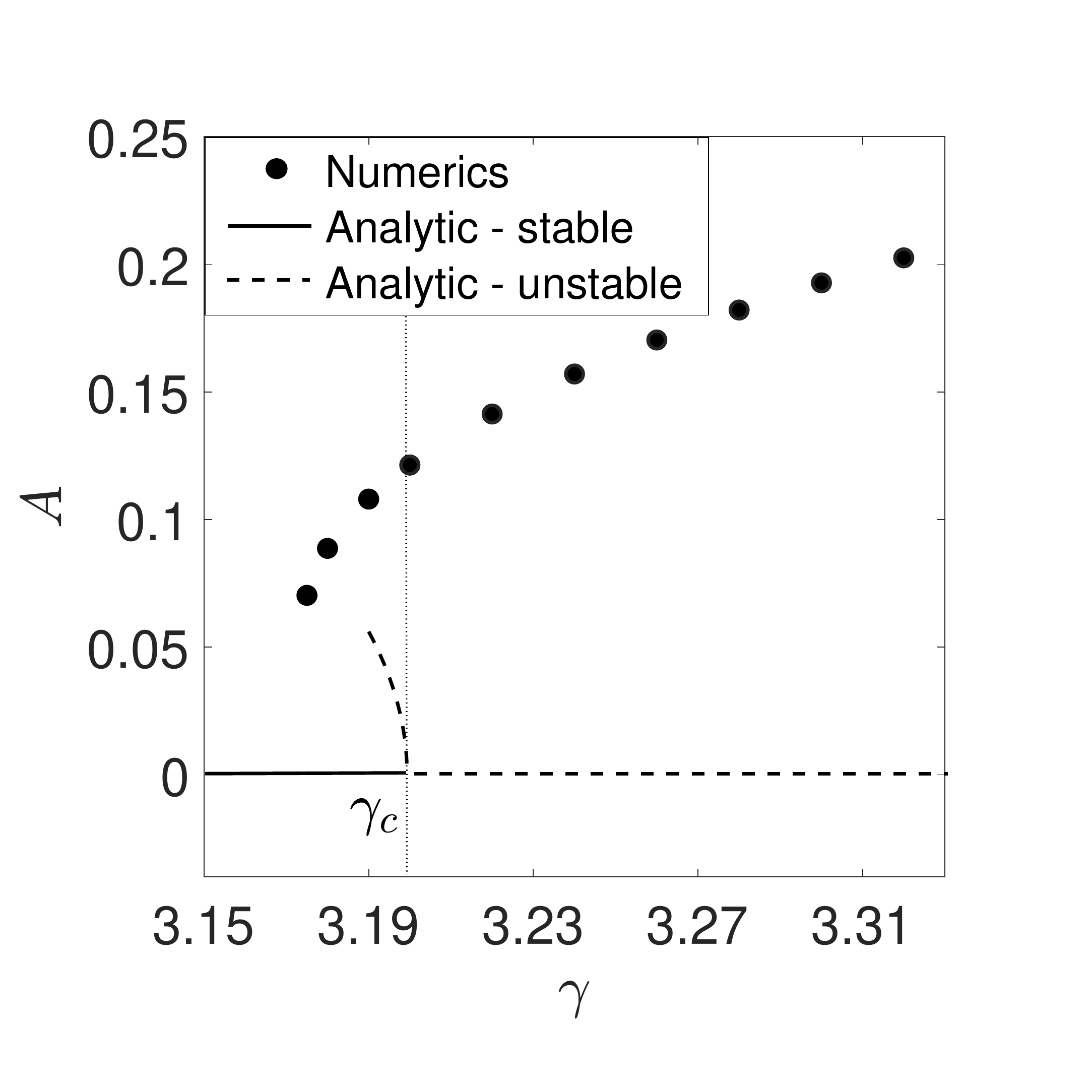}}
	\subfloat[ $ L=5 $]
	{\includegraphics[width=.33\textwidth]{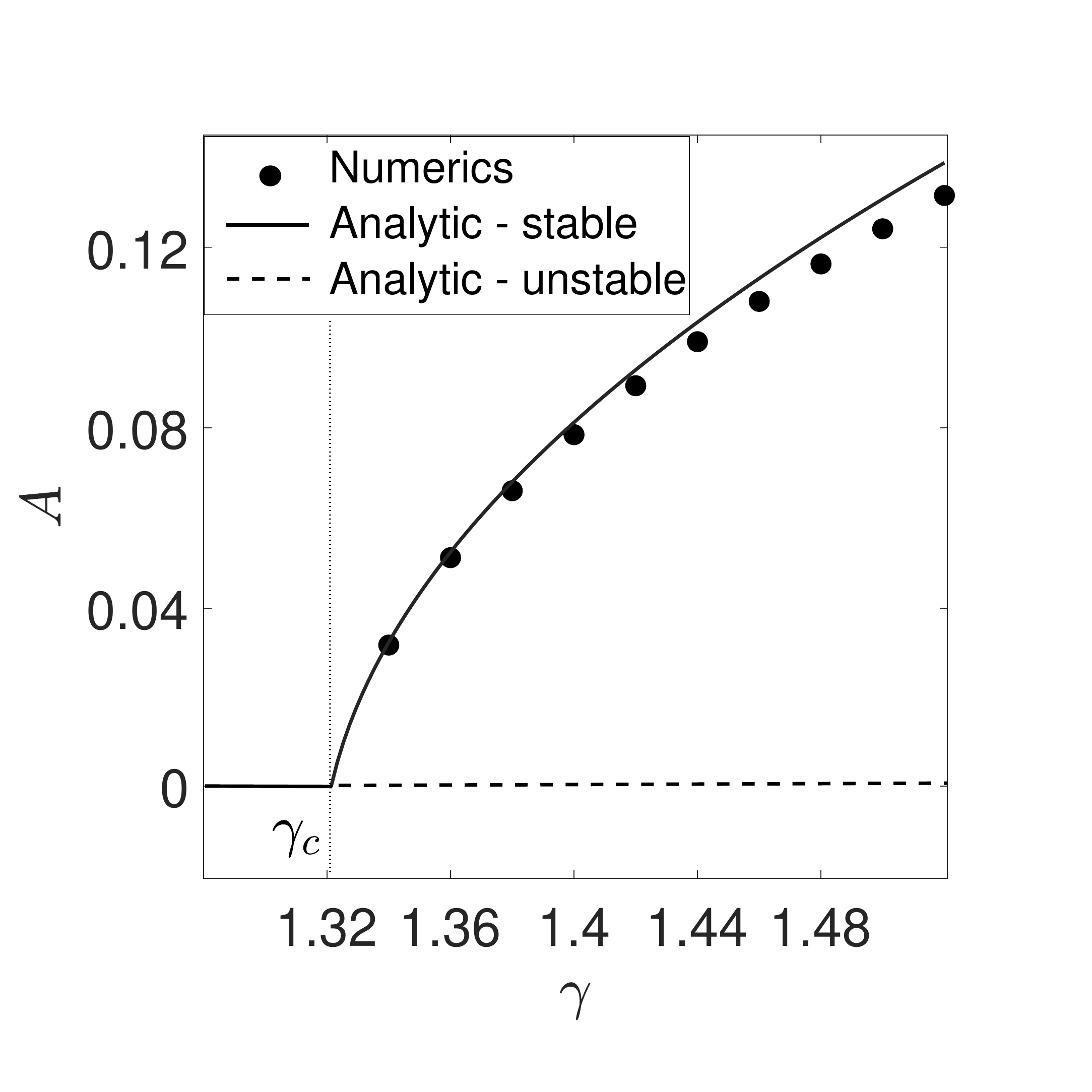}}
	\subfloat[$ L=15 $]
	{\includegraphics[width=.33\textwidth]{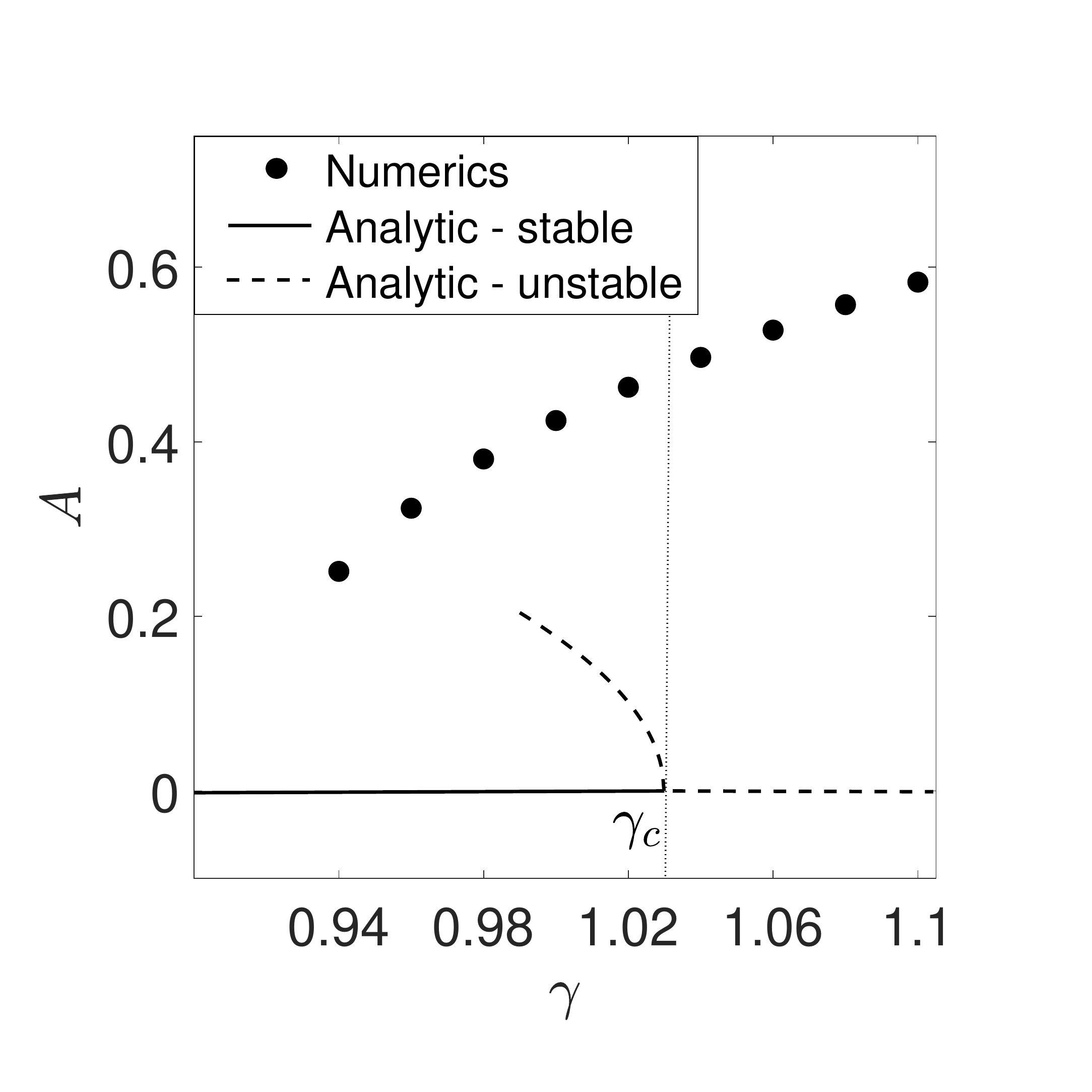}}
		\caption{Comparison between analytical and numerical bifurcation diagrams of system \eqref{eq:ndsystem} with density-dependent advection strength $\gamma>0$, and nonlocal kernel $K=K_1$ (see Equation \eqref{eq:ndtop-hat}), $\bar{u}_1=0.1$ and $\bar{u}_2=10$, for different values of the length of the domain $L$. These scenarios correspond to Figure \ref{fig:L} (b) (center). Dashed and solid lines represent unstable and stable branches, respectively, which are computed analytically, while the dots are computed numerically. As the length of the domain increases, the system changes its qualitative behaviour.  In (a): $L=2.7$ and the system exhibits a subcritical bifurcation at $\gamma=\gamma_c=3.19933$. In (b), $L=5$ and at $\gamma=\gamma_c=1.32131$ a branch of stable solutions bifurcates from the homogeneous state . In (c), $L=15$ and the system exhibits a subcritical bifurcation at $\gamma=\gamma_c=1.02985$.}
	\label{fig:avoid_sub_super_sub}
    \end{figure}

Figure \ref{fig:avoid_sub_super_sub} shows bifurcation diagrams obtained by fixing $\bar{u}_1=0.1$ and $\bar{u}_2=10$ and by changing $L$, in the mutual avoidance regime ($\gamma>0$). This case corresponds to the scenario shown in Figure \ref{fig:L}(b) (center). Dashed and solid lines represent unstable and stable branches, respectively, computed analytically, while the dots are computed numerically. For $L=2.7$, the weakly nonlinear analysis predicts a subcritical bifurcation, and the numerical simulations confirm this result. In fact, just beyond the instability threshold ($\gamma>\gamma_c\approx3.20$), we find stable large amplitude solutions, which persist when we decrease the control parameter $\gamma$ below the instability threshold (Figure \ref{fig:avoid_sub_super_sub}(a)). For $L=5$, the analysis predicts a supercritical bifurcation and, again, the numerical simulations confirm this result. In Figure \ref{fig:avoid_sub_super_sub}(b) we see, indeed, a good matching between the analytical branch and the numerical solutions, as long as $\gamma$ is sufficiently close to the bifurcation threshold $\gamma_c\approx1.32$. Finally, for $L=15$ the subcritical bifurcation predicted by our analysis is also detected numerically (see Figure \ref{fig:avoid_sub_super_sub}(c)).  Here, we observe  bistability between the homogeneous steady state and non-homogeneous solutions below the instability threshold $\gamma_c\approx1.03$.

 
\begin{figure}[H]
		\centering
        $\gamma<0$, $\bar{u}_1\neq \bar{u}_2$
          \subfloat[$L=2.5$]
        {\includegraphics[width=0.33\textwidth]{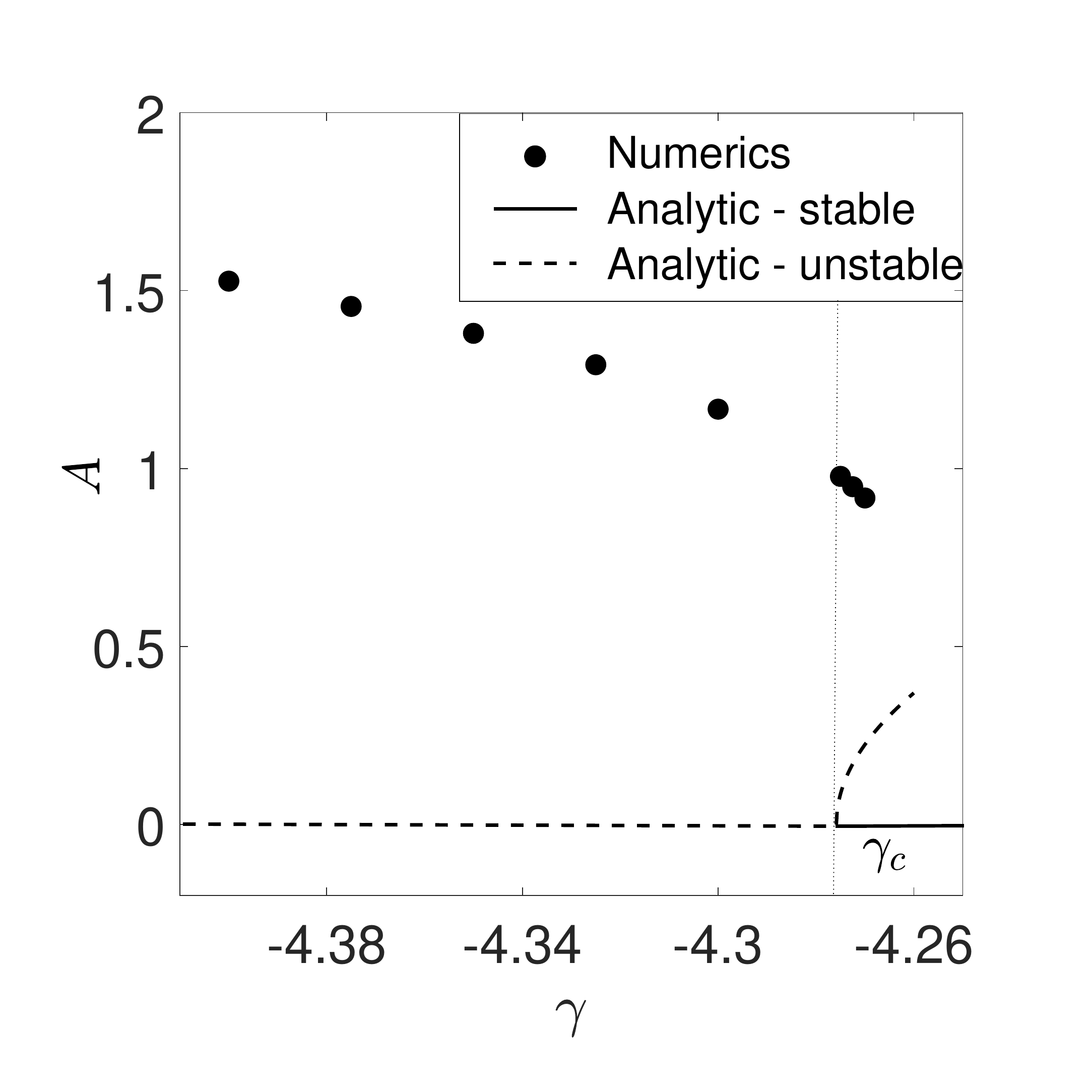}}
	\subfloat[ $ L=5 $]
	{\includegraphics[width=.33\textwidth]{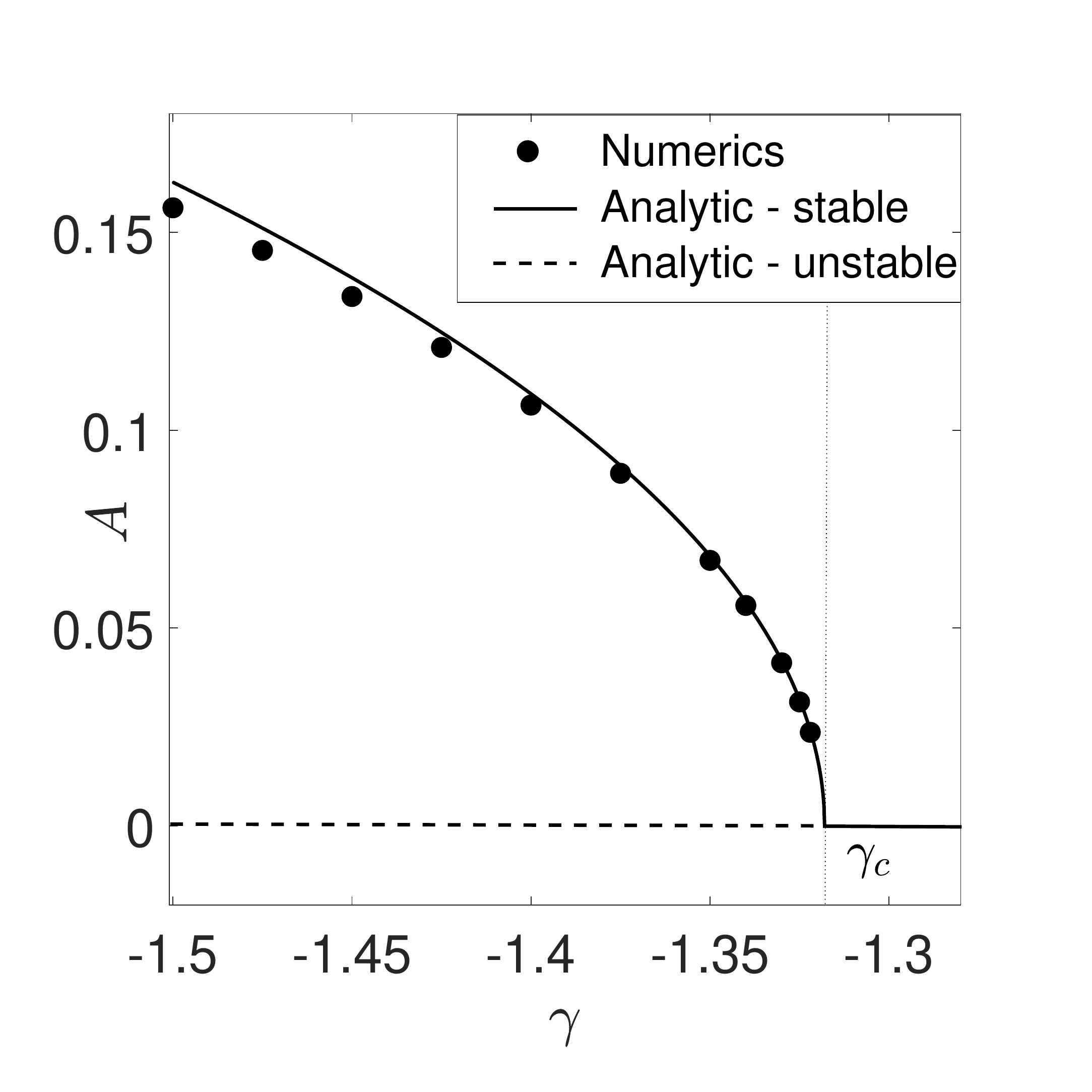}}	
         \subfloat[ $ L=10 $]
	{\includegraphics[width=.33\textwidth]{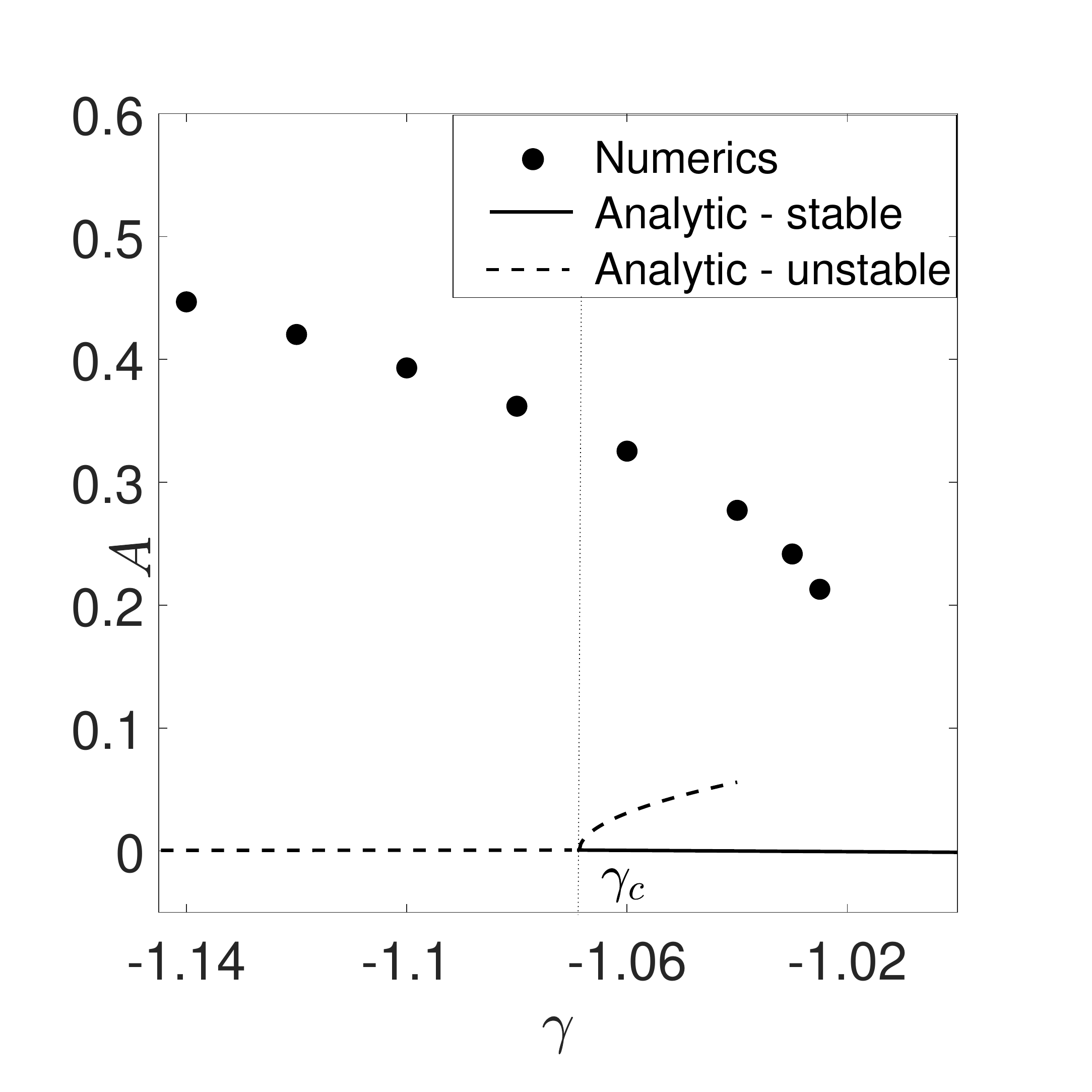}}
		\caption{Comparison between analytical and numerical bifurcation diagrams of system \eqref{eq:ndsystem} with $\gamma<0$, $K=K_1$ (see Equation \eqref{eq:ndtop-hat}), $\bar{u}_1=0.1$ and $\bar{u}_2=10$, for different values of the length of the domain $L$. These scenarios correspond to Figure \ref{fig:L} (d) (center). Dashed and solid lines represent unstable and stable branches, respectively, which are computed analytically, while the dots are computed numerically. As the length of the domain increases, the system changes its qualitative behaviour.  In (a): $L=2.5$ and the system exhibits a subcritical bifurcation at $\gamma=\gamma_c=-4.2758$. In (b), $L=5$ and at $\gamma=\gamma_c=-1.32131$ a branch of stable solutions bifurcates from the homogeneous state. In (c), $L=10$ and the system exhibits a subcritical bifurcation at $\gamma=\gamma_c=-1.06895$.}
	\label{fig:attr_sub_super_sub}
    \end{figure}

Figure \ref{fig:attr_sub_super_sub} shows bifurcation diagrams obtained by fixing $\bar{u}_1=0.1$ and $\bar{u}_2=10$, for three different values of $L$, in the mutual avoidance regime ($\gamma<0$). This case corresponds to the scenario shown in Figure \ref{fig:L}(d) (center). The numerical simulations, again, confirm the results of the weakly nonlinear analysis: we have detected subcritical transitions for $L=2$ and $L=10$, and a stable branch bifurcating supercritically for $L=5,$ whose amplitude is well approximated by the weakly nonlinear analysis.
    
    It remains to analyze the case $\bar{u}_1=\bar{u}_2$, corresponding to the scenarios depicted in Figures \ref{fig:L}(a) ($\gamma>0$) and (c) ($\gamma<0$). In this case, three different types of bifurcation are predicted by the analysis: subcritical bifurcations (for $\Lambda<0$), unstable supercritical bifurcations (for $\Lambda>0$ and $\Gamma<0$) and stable supercritical bifurcations (for $\Lambda>0$ and $\Gamma>0$) (see Figure \ref{fig:gammac_vs_L}). In particular, for $\gamma>0$, system \eqref{eq:ndsystem} undergoes subcritical bifurcations for $2<L<3$, unstable supercritical bifurcations for $3<L<3.5$, and stable supercritical bifurcations for $L>3.5$ (see Figure \ref{fig:gammac_vs_L}(a)). 
    
\begin{figure}[H]
        \centering
        $\gamma>0$, $\bar{u}_1 = \bar{u}_2$
		\subfloat[$L=3.1$]
		{\includegraphics[width=1\textwidth]{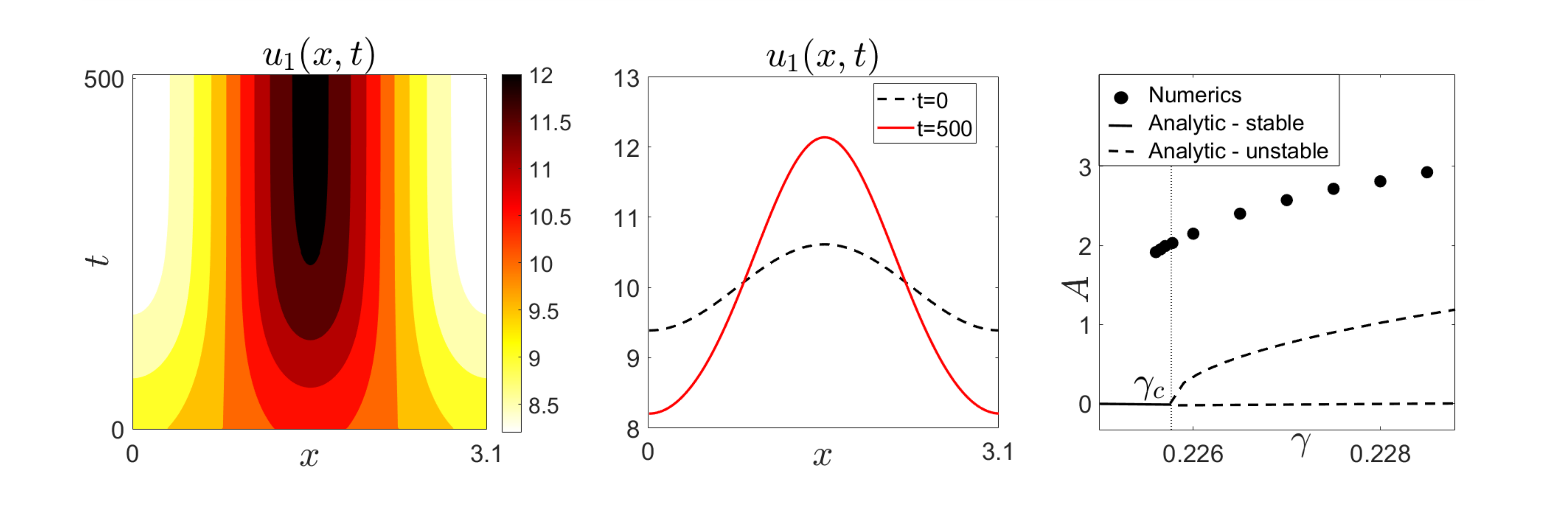}}\\
		\subfloat[$L=4$]
		{\includegraphics[width=1\textwidth]{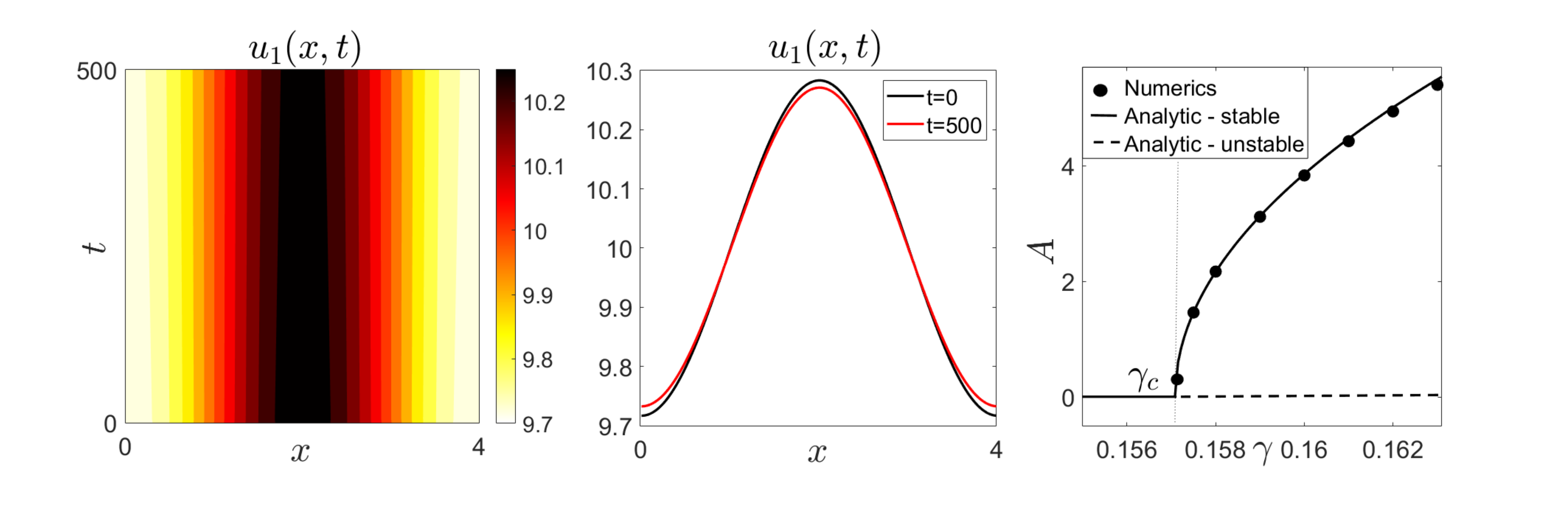}}
		\caption{Numerical investigation of system \eqref{eq:ndsystem} in the mutual avoidance regime ($ \gamma>0 $) with $\bar u_1=\bar u_2=10$, for two different values of $L$. In (a): $L=3.1$ and the analysis predict an unstable supercritical bifurcation at $\gamma=\gamma_c=0.225754$. On the left, numerical simulation showing that the system moves away from the unstable solution and evolves toward a large amplitude pattern. In the center, initial condition and the final stationary state. On the right, comparison between analytical and numerical bifurcation diagrams.
        In (b): $L=4$ and the analysis predict a stable supercritical bifurcation at $\gamma=\gamma_c=0.15708$. On the left, numerical simulation showing that the system moves towards the stable small amplitude solution. In the center, initial condition and the final stationary state. On the right, comparison between analytical and numerical bifurcation diagrams.}
		\label{fig:Stable_Unstable_branches}
	\end{figure}

    In Figure \ref{fig:Stable_Unstable_branches} we analyze System \eqref{eq:ndsystem} with $\gamma>0$ and $\bar u_1=\bar u_2=10$, for $L=3.1$ in (a), and $L=4$ in (b). In Figure \ref{fig:Stable_Unstable_branches}(a) (left) we show the spatio-temporal evolution of a numerical solution whose initial condition is a small perturbation of the weakly nonlinear solution with $L=3.1$. We observe that the numerical solution moves away from the initial condition and evolves toward a large amplitude pattern. The initial condition and the final stationary state are reported in Figure \ref{fig:Stable_Unstable_branches}(a) (center). Therefore, when the supercritical branch is unstable, the system supports large amplitude patterns. These solutions exist even below the bifurcation threshold, as shown by the bifurcation diagram in Figure \ref{fig:Stable_Unstable_branches}(a) (right). These large amplitude solutions are not predicted by the weakly nonlinear analysis.  However we conjecture that they might be obtained analytically by expanding the weakly nonlinear analysis to higher orders. 
    
    In Figure \ref{fig:Stable_Unstable_branches}(b) (left) we show the spatio-temporal evolution of a numerical solution whose initial condition is a small perturbation of the weakly nonlinear solution with $L=4$. In this case, the analysis predicts that the small amplitude pattern is stable. In the numerical simulation we observe that the solution moves towards a small amplitude pattern, which is well approximated by the weakly nonlinear analysis. This result confirms the stability predicted by our analysis. The initial condition and the final stationary state are reported in Figure \ref{fig:Stable_Unstable_branches}(b) (center). Finally, a comparison between the analytical and numerical bifurcation diagrams is shown in \ref{fig:Stable_Unstable_branches}(b) (right).

\begin{figure}[H] 
  \centering
        $\gamma<0$, $\bar{u}_1 = \bar{u}_2$\\
          \subfloat[$L=5$]
        {\includegraphics[width=0.33\textwidth]{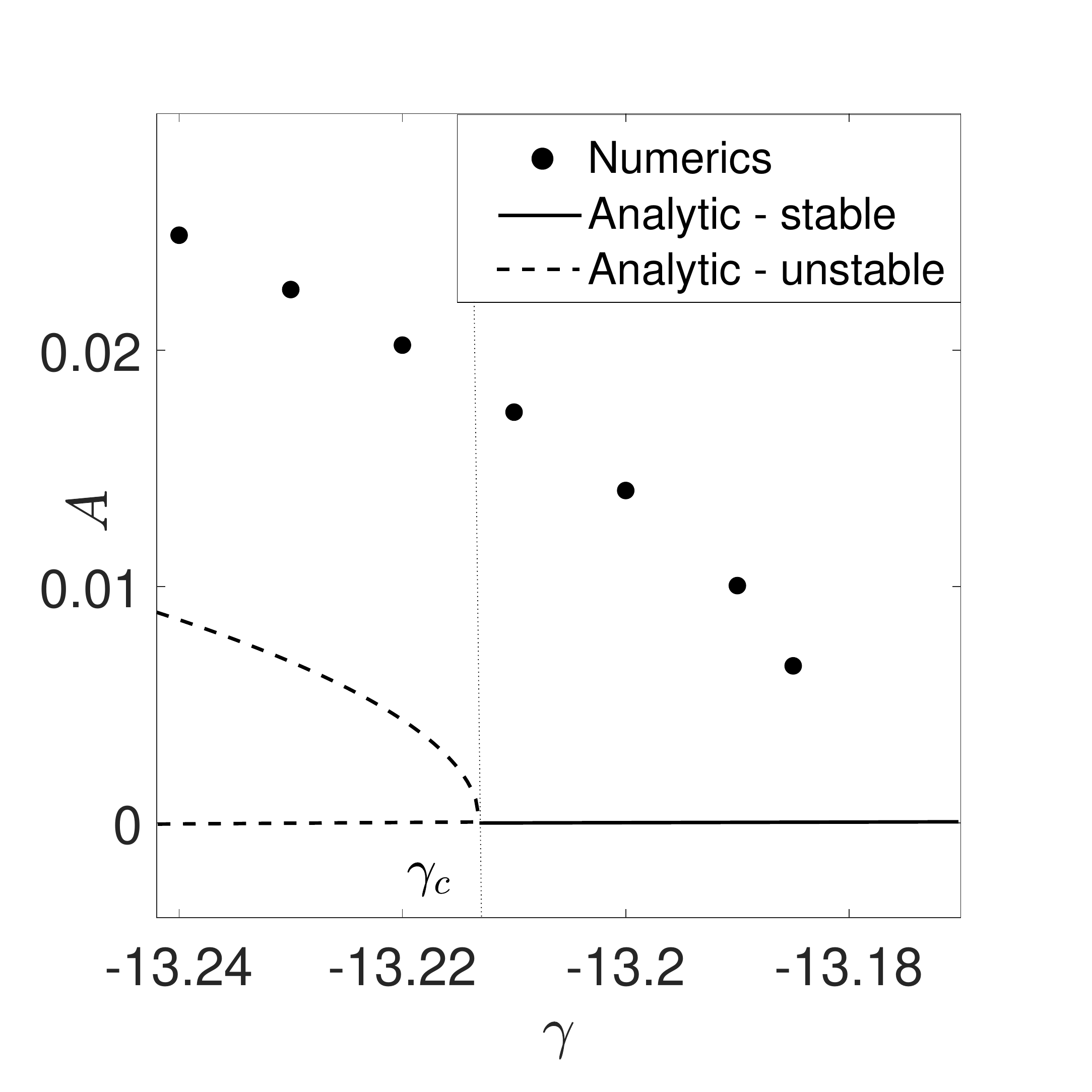}}
	\subfloat[ $ L=10 $]
	{\includegraphics[width=.33\textwidth]{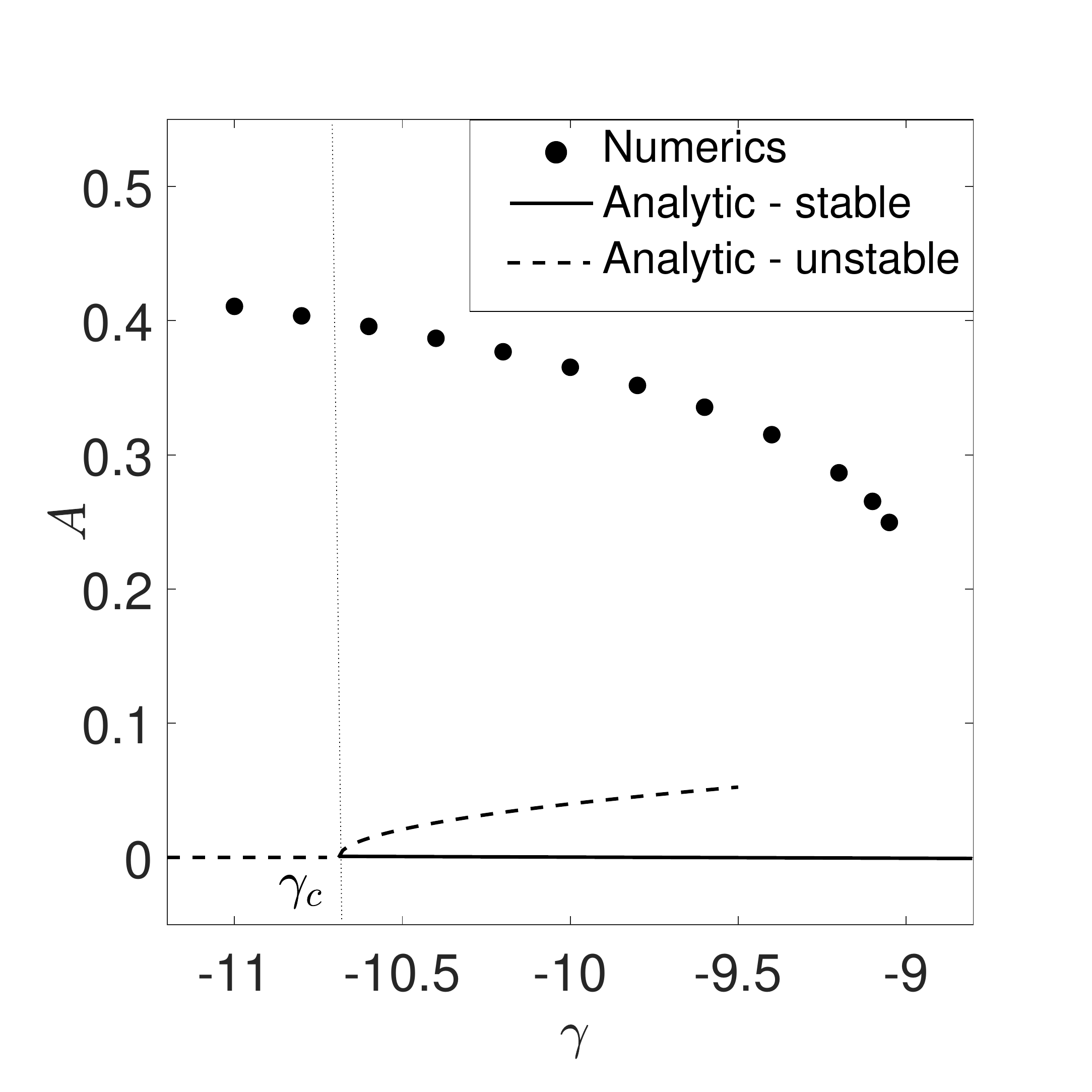}}	
		\caption{Comparison between analytical and numerical bifurcation diagrams of system \eqref{eq:ndsystem} with $\gamma<0$, $K=K_1$ (see Equation \eqref{eq:ndtop-hat}), $\bar{u}_1=\bar{u}_2=10$, for different values of the length of the domain $L$. These scenarios correspond to Figure \ref{fig:L} (c) (right). Dashed and solid lines represent unstable and stable branches, respectively, which are computed analytically, while the dots are computed numerically. As the length of the domain increases, the system changes its qualitative behaviour.  In (a): $L=5$ and the system exhibits a supercritical bifurcation at $\gamma=\gamma_c\approx-13.2$, giving rise to a branch of unstable small amplitude solutions. In (b), $L=10$ and at $\gamma=\gamma_c\approx-10.7$ the system exhibits a subcritical bifurcation}
	\label{fig:super_sub_attr}
    \end{figure}

       Finally, Figure \ref{fig:super_sub_attr} shows analytic and numerical bifurcation diagrams of System \eqref{eq:ndsystem} with $\gamma<0$ and $\bar u_1=\bar u_2=10$. Our previous analysis predicts unstable supercritical bifurcations for $2<L<6$, and subcritical bifurcations for $L>6$ (see Figure \ref{fig:gammac_vs_L} (b)). We have verified these results numerically, and the comparisons between analytical and numerical bifurcation diagrams are shown in Figure \ref{fig:super_sub_attr}.

\subsection{Bistability between small amplitude patterns and strongly modulated solutions}\label{sec:th_bistability}

    The existence of non-constant solutions to system \eqref{eq:ndsystem}, far away from any bifurcation of the constant solution, was already detected and analyzed in \cite{GHLP22} {using an energy method}. By minimising an energy functional associated with the system, nontrivial stationary solutions were revealed which, as $L$ increases, tend to look increasingly like piecewise constant functions, when $\gamma>0$, or spike solutions, when $\gamma<0$.  We call such solutions \textit{strongly modulated} because they are given by the superposition of more than one unstable Fourier mode. In this section, we will combine numerical and analytic solutions inferred by both the weakly nonlinear analysis here and the results presented in \cite{GHLP22} to construct more comprehensive bifurcation diagrams.

For this, we focus on the case $\gamma>0$ and $\bar{u}_1=\bar{u}_2$. Here, the system exhibits supercritical bifurcations for large values of $L$ (see Figure \ref{fig:L} (a)). Also, as shown in Figure \ref{fig:Gamma} (a), these supercritical bifurcations generate stable small amplitude patterns. In \cite{GHLP22} we showed that under the same conditions (that is $L\gg 1$, $\gamma>0$ and $\bar{u}_1=\bar{u}_2$), the system supports strongly modulated patterns.  
Therefore we expect that for $L$ sufficiently large, there exist parameter regions in which small amplitude patterns and strongly modulated solutions coexist and are stable. 

\begin{figure}[H]
		\centering
        \subfloat[$L=10$]
        {\includegraphics[scale=0.33]{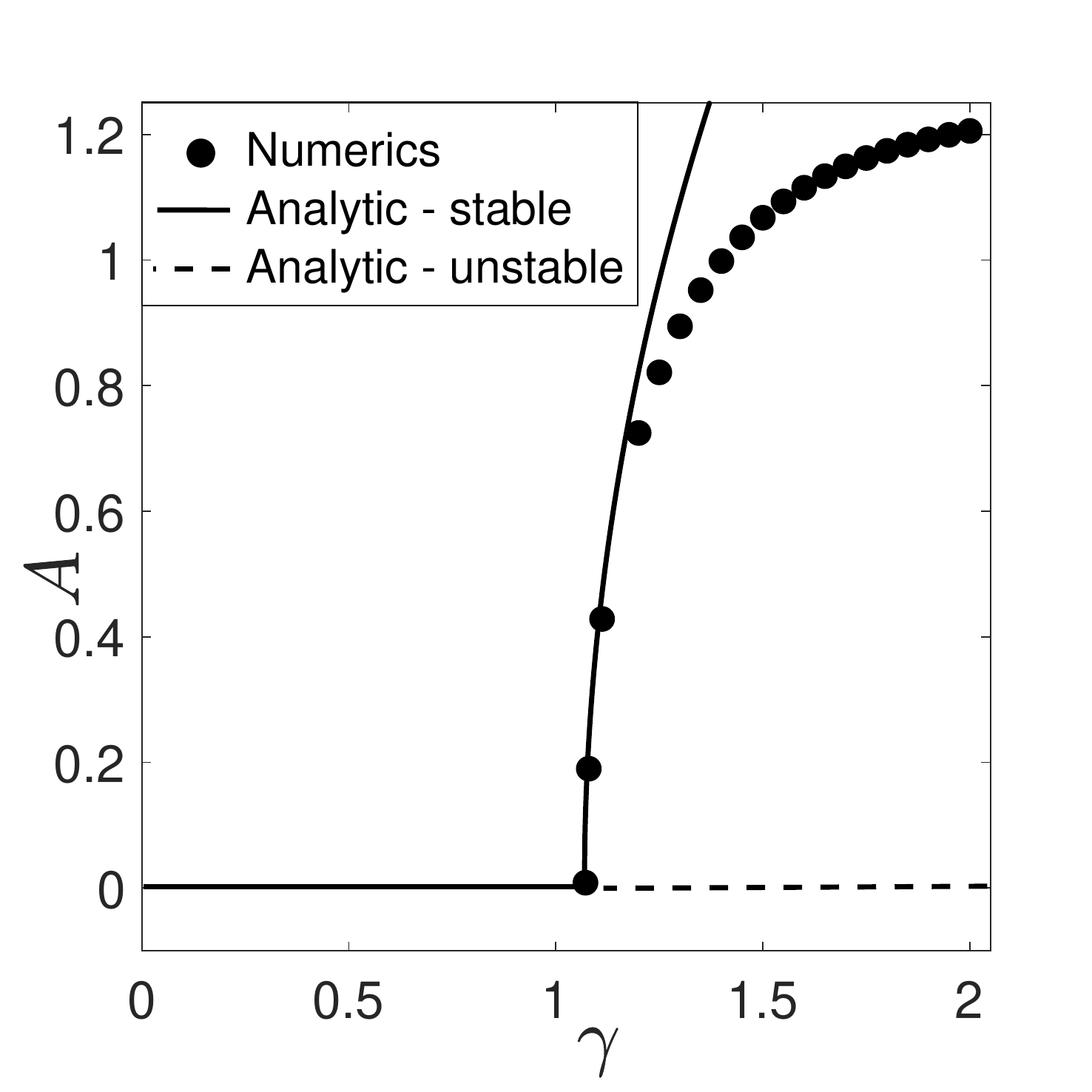}}
          \subfloat[$L=20$]
        {\includegraphics[scale=0.33]{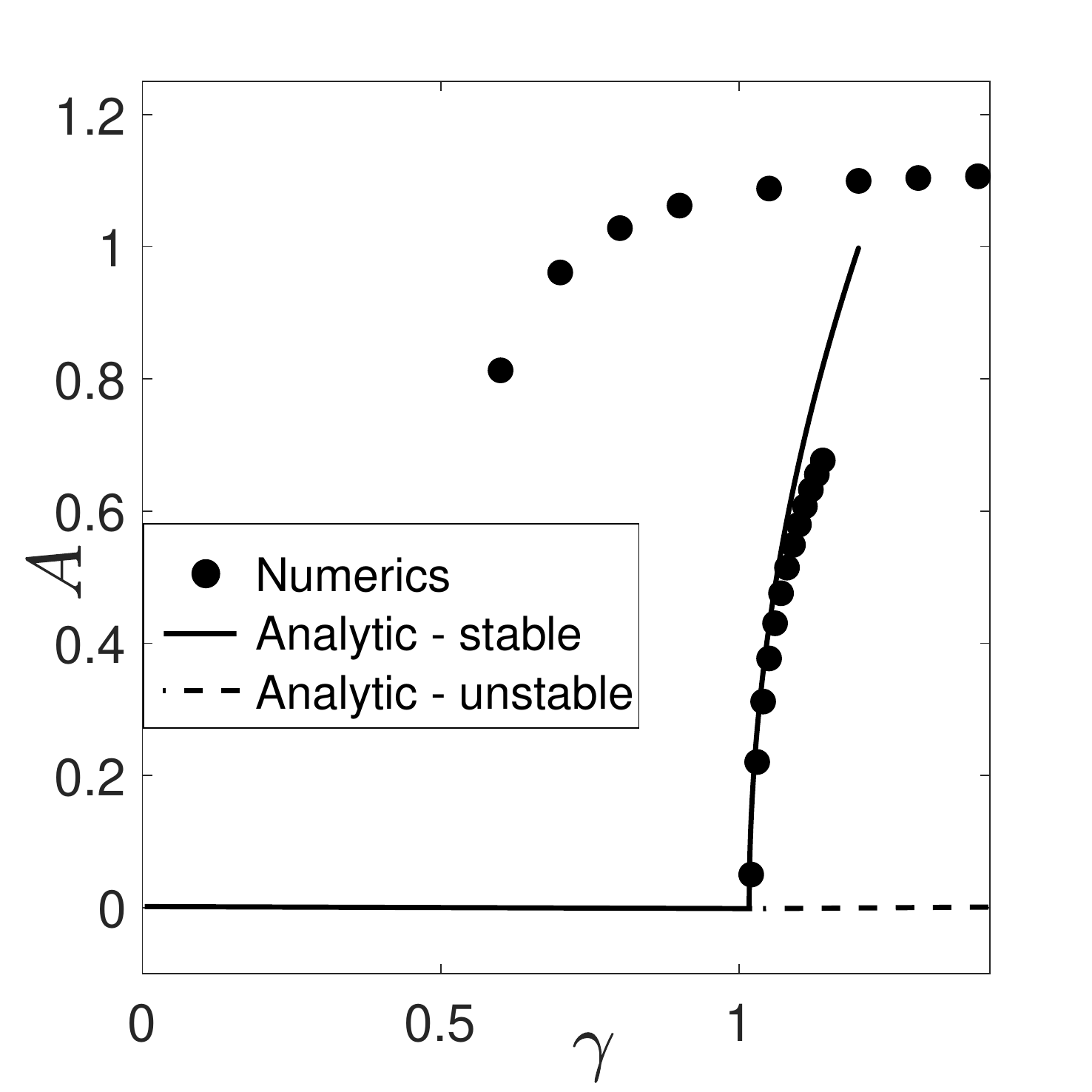}}
	\subfloat[ $ L=50 $]
	{\includegraphics[scale=0.33]{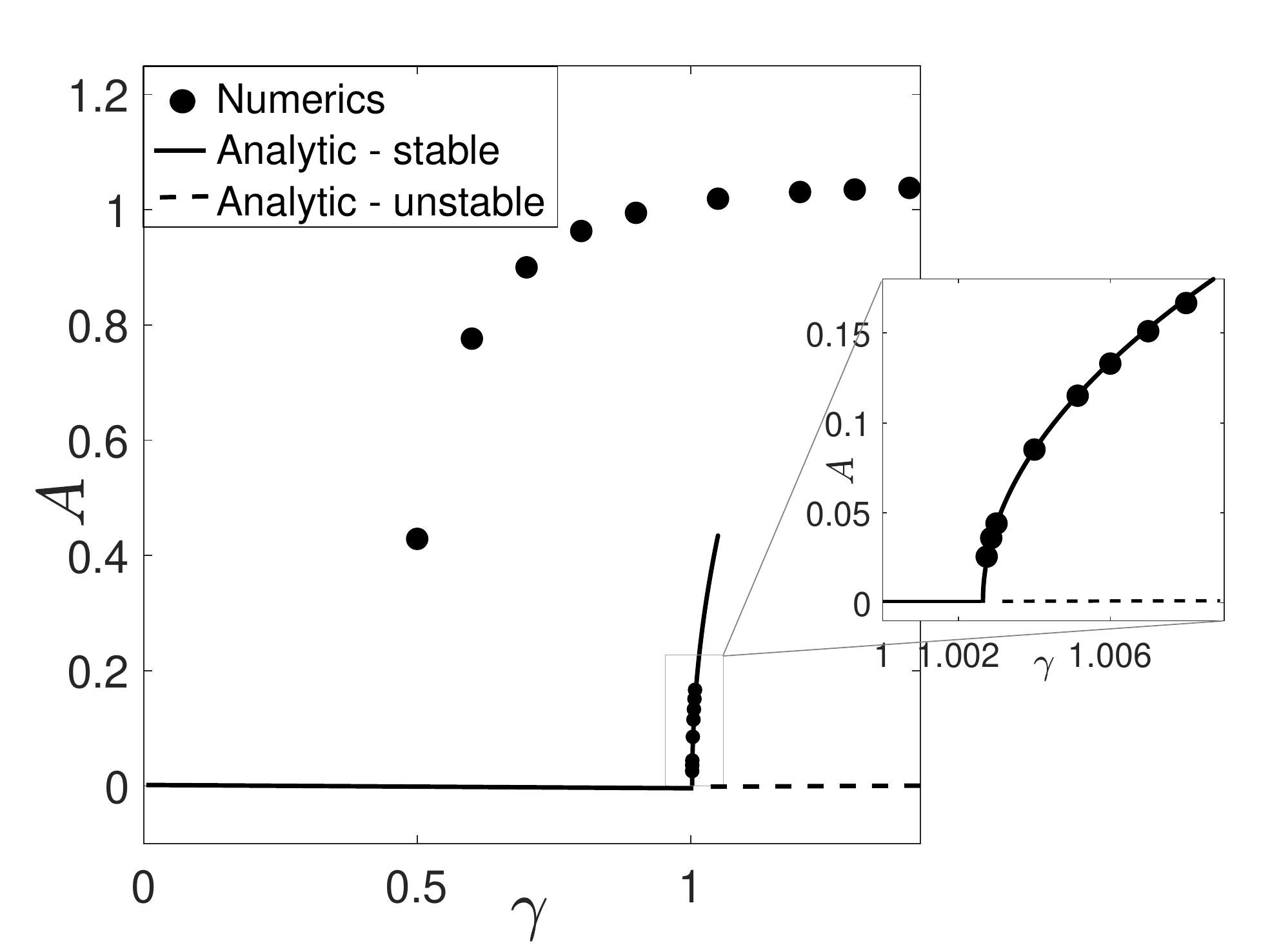}}	
		\caption{Bifurcation diagrams of system \eqref{eq:ndsystem} with $\gamma>0$, $K=K_1$ (see Equation \eqref{eq:ndtop-hat}), $\bar{u}_1=\bar{u}_2=1$, for different values of the length of the domain $L$. These scenarios correspond to Figure \ref{fig:L} (a) (center). Dashed and solid lines represent unstable and stable branches, respectively, which are computed analytically, while the dots are computed numerically.  The system exhibits a supercritical stable bifurcation at: $\gamma=\gamma_c=1.06896$ in (a);  $\gamma=\gamma_c=1.01664$ in (b);  $\gamma=\gamma_c=1.0264$ in (c). As $L$ becomes sufficiently large, the system support strongly modulated patterns which coexist with stable small amplitude patterns.}
	\label{fig:Bistability}
    \end{figure}

We have verified this numerically and the results are shown in Figure \ref{fig:Bistability}. When $L$ is not too large, the system admits small amplitude solutions that bifurcate supercritically from the homogeneous steady state and remains stable as the control parameter $\gamma$ increases (see Figure \ref{fig:Bistability} (a)). In this case, we do not find strongly modulated solutions. As $L$ increases, the supercritical branch of patterns predicted by the weakly nonlinear analysis still exists and is stable as long as $\gamma$ is sufficiently close to the bifurcation threshold (see Figure \ref{fig:Bistability} (b)). However, a second branch appears higher up, representing the strongly modulated solutions predicted by \cite{GHLP22}. 
As $L$ increases further, the branch of stable small amplitude solutions becomes smaller and smaller (Figure \ref{fig:Bistability}(c)), and the solutions transition to strongly modulated for values of $\gamma$ closer to the bifurcation threshold.

    \section{Discussion}\label{sec:discussion}
    We have analysed bifurcations for a nonlocal advection diffusion system with two interacting populations that either mutually avoid or mutually attract. First, we analysed the linear stability of the homogeneous steady state and recovered the instability thresholds. Beyond these thresholds, the homogeneous steady state becomes unstable and the system is expected to form spatially inhomogeneous patterns. To predict the evolution of the system in the unstable regime, we used weakly nonlinear analysis to recover the equations governing the amplitude of the pattern and approximations of the inhomogeneous solutions. We  found that the amplitude equations consist of a Ginzburg-Landau equation coupled with an equation for the zero mode. Indeed, we  obtained a sequence of linear problems whose general solutions must be a linear combination of the critical mode and the zero mode. This follows from the fact that the system under study obeys a conservation law. An equivalent result was shown in \cite{Murray01}, where similar amplitude equations were derived using symmetry and scaling arguments. By means of the amplitude equations, we recovered the condition that ensures the stability of the patterns bifurcating from the homogeneous steady state. 
    
    To obtain concrete numerical results, we analysed the case where the spatial-averaging kernel, $K$, is a top-hat distribution. By combining analysis of the amplitude equation with numerical solutions, we showed that the system exhibits a variety of different types of bifurcations and bistability regimes, strongly depending on the ratio $l/\alpha$. In particular, we found stable small amplitude patterns bifurcating supercritically from the homogeneous steady state at the onset of the instability. We also found subcritical regimes generating unstable small amplitude patterns, which coexist with both the stable homogeneous solution and stable large amplitude patterns. In this case, numerics revealed an hysteresis effect due to the bistability between two stationary states. Finally, we also found supercritical bifurcations generating unstable small amplitude patterns. Beyond the instability threshold, we numerically detected  stable large amplitude patterns that persist even when decreasing the bifurcation parameter below the instability threshold, revealing again a hysteresis effect similar to that found in the subcritical regime. 
    
    By combining weakly nonlinear analysis, numerical simulations and the energy functional analysis from  \cite{GHLP22}, {we obtained a comprehensive bifurcation picture. We found} parameter regions exhibiting bistability between small amplitude patterns and strongly modulated solutions, when $l/\alpha\gg 1$. The range of bistability becomes smaller and smaller as $l/\alpha$ increases, because the small amplitude patterns lose their stability for values of the control parameter increasingly closer to the bifurcation threshold (Figure \ref{fig:Bistability}).  Overall, our analysis reveals that our system may display discontinuous phase transitions either when $\alpha \approx l$ or when the sensing range $\alpha$ is very small compared to the length $l$ of the domain. 
    
    Our study provides an example of how to combine different and complementary approaches to recover more comprehensive pictures of the bifurcation diagrams. To extend these results further, it would be interesting to expand the weakly nonlinear analysis up to higher orders.  Such an approach could reveal analytically some of the large amplitude branches here found numerically, as well as the branches of solutions connecting small and large amplitude patterns. Numerical continuation software, such as pde2path \cite{uecker2021numerical}, gives another way of approaching this problem \cite{soresina2022hopf,consolo2022eckhaus}. Our analysis revealed parameter regions with bistability between two extended states, a scenario in which systems often exhibit snaking branches of localized solutions \cite{burke2007snakes, verschueren2017model}. Extending our weakly nonlinear analysis to higher orders may help locate the codimension-two point where the nascence of localised structures may take place, which would be an interesting subject for future work.

    Our focus here has been on a particular example of Equation \eqref{eq:system} \cite{PL19}, with just two populations and no self-interaction terms ($N=2$, $\gamma_{ii}=0$). However, even in this relatively-simple system, we found an unexpectedly rich variety of patterning scenarios. Therefore, we conjecture that  analysis of the system with $N \geq 3$ populations and/or $\gamma_{ii}\neq 0$ would reveal even more complex patterning and bifurcation structure. Our next goal, indeed, is to analyse the more general scenarios ($N \geq 3$, $\gamma_{ii}\neq 0$). 
    A possible way forward might be to analyse phase transitions by combining the tools used here with those from \cite{carrillo2020long}. In \cite{carrillo2020long} the authors studied the phase transitions of the Mckean-Vlasov equation by analysing the minimizers of the energy associated to the problem. Combining this with weakly nonlinear analysis might shed light on the number of steady states at the onset of an instability, and consequently on type of phase transition occurring when the bifurcation parameter crosses the instability threshold.  

    System \eqref{eq:system} has several applications to natural systems and, in particular, to ecological systems. Therefore the analysis presented in this paper, as well as possible future extensions, might help to address some important ecological questions regarding the emergence of territories, as well as their sizes and stability  \cite{potts2014animal}. Indeed, variations in territory size and shape can strongly affect population structure and dynamics \cite{adams2001approaches}, therefore understanding the mechanisms and consequences of these changes is crucial for informing the design of efficient conservation strategies. Our results support the hypothesis that the formation of territorial patterns is not just a consequence of a heterogeneity in resources distribution, but that they can emerge as a consequence of animal behaviour and mutual interactions \cite{adams2001approaches,giuggioli2011animal, potts2014animal}. Our analysis also predicts that a small sensing range {relative to the length of the domain} can facilitate a territory instability, in agreement with other theoretical studies suggesting that poor sensory information can promote the range size instability (\cite{tao2016dynamic}).
   In summary, the analysis of the class of models \eqref{eq:system} with the techniques here presented and discussed can help to resolve biological and ecological questions that may be inaccessible to experimental investigation.
   
   \bigskip
   \noindent{\bf Acknowledgements:} JRP and VG acknowledge support of Engineering and Physical Sciences Research Council (EPSRC) grant EP/V002988/1 awarded to JRP. VG is also grateful for support from the National Group of Mathematical Physics (GNFM-INdAM).
   TH is supported through a discovery grant of the Natural Science and Engineering Research Council of Canada (NSERC), RGPIN-2017-04158. MAL gratefully acknowledges support from NSERC Discovery Grant RGPIN-2018-05210 and from the Gilbert and Betty Kennedy Chair in Mathematical Biology.

   \noindent{\bf Declarations of interest:} The authors have no competing interests to declare.
   



\bibliographystyle{plain} 
 \bibliography{bib_wnl}

\begin{thebibliography}{10}

\bibitem{adams2001approaches}
Eldridge~S Adams.
\newblock Approaches to the study of territory size and shape.
\newblock {\em Annual Review of Ecology and Systematics}, 32(1):277--303, 2001.

\bibitem{basaran2022large}
Mustafa Basaran, Y~Ilker Yaman, Tevfik~Can Y{\"u}ce, Roman Vetter, and Askin
  Kocabas.
\newblock Large-scale orientational order in bacterial colonies during inward
  growth.
\newblock {\em Elife}, 11:e72187, 2022.

\bibitem{bergeon2008eckhaus}
Alain Bergeon, J~Burke, E~Knobloch, and I~Mercader.
\newblock Eckhaus instability and homoclinic snaking.
\newblock {\em Physical Review E}, 78(4):046201, 2008.

\bibitem{bilotta2018eckhaus}
Eleonora Bilotta, Francesco Gargano, Valeria Giunta, Maria~Carmela Lombardo,
  Pietro Pantano, and Marco Sammartino.
\newblock Eckhaus and zigzag instability in a chemotaxis model of multiple
  sclerosis.
\newblock {\em Atti della Accademia Peloritana dei Pericolanti-Classe di
  Scienze Fisiche, Matematiche e Naturali}, 96(S3):9, 2018.

\bibitem{BLP02}
Brian~K Briscoe, Mark~A Lewis, and Stephen~E Parrish.
\newblock Home range formation in wolves due to scent marking.
\newblock {\em Bulletin of Mathematical Biology}, 64(2):261--284, 2002.

\bibitem{Stevens}
M.~Burger, M.~Di~Francesco, S.~Fagioli, and A.~Stevens.
\newblock Sorting phenomena in a mathematical model for two mutually
  attracting/repelling species.
\newblock {\em SIAM Journal on Mathematical Analysis}, 50(3):3210--3250, 2018.

\bibitem{burke2007snakes}
John Burke and Edgar Knobloch.
\newblock Snakes and ladders: localized states in the {S}wift--{H}ohenberg
  equation.
\newblock {\em Physics Letters A}, 360(6):681--688, 2007.

\bibitem{buttenschon2021non}
Andreas Buttensch{\"o}n and Thomas Hillen.
\newblock {\em Non-local cell adhesion models}.
\newblock Springer, 2021.

\bibitem{carrillo2020long}
JA~Carrillo, RS~Gvalani, GA~Pavliotis, and A~Schlichting.
\newblock Long-time behaviour and phase transitions for the mckean--vlasov
  equation on the torus.
\newblock {\em Archive for Rational Mechanics and Analysis}, 235(1):635--690,
  2020.

\bibitem{CKY19}
Jos{\'e}~A Carrillo, Katy Craig, and Yao Yao.
\newblock Aggregation-diffusion equations: dynamics, asymptotics, and singular
  limits.
\newblock {\em Active Particles, Volume 2: Advances in Theory, Models, and
  Applications}, pages 65--108, 2019.

\bibitem{consolo2022eckhaus}
Giancarlo Consolo and Gabriele Grif{\'o}.
\newblock Eckhaus instability of stationary patterns in hyperbolic
  reaction--diffusion models on large finite domains.
\newblock {\em Partial Differential Equations and Applications}, 3(5):57, 2022.

\bibitem{CH93}
Mark~C Cross and Pierre~C Hohenberg.
\newblock Pattern formation outside of equilibrium.
\newblock {\em Reviews of Modern Physics}, 65(3):851, 1993.

\bibitem{CG09}
Michael Cross and Henry Greenside.
\newblock {\em Pattern formation and dynamics in nonequilibrium systems}.
\newblock Cambridge University Press, 2009.

\bibitem{eisenbach2004chemotaxis}
Michael Eisenbach, A~Tamada, GM~Omann, JE~Segall, RA~Firtel, R~Meili, David
  Gutnick, Mazal Varon, Joseph~W Lengeler, and Fujio Murakami.
\newblock {\em Chemotaxis}.
\newblock World Scientific Publishing Company, 2004.

\bibitem{ellison2020mechanistic}
Natasha Ellison, Ben~J Hatchwell, Sarah~J Biddiscombe, Clare~J Napper, and
  Jonathan~R Potts.
\newblock Mechanistic home range analysis reveals drivers of space use patterns
  for a non-territorial passerine.
\newblock {\em Journal of Animal Ecology}, 89(12):2763--2776, 2020.

\bibitem{giuggioli2011animal}
Luca Giuggioli, Jonathan~R Potts, and Stephen Harris.
\newblock Animal interactions and the emergence of territoriality.
\newblock {\em PLoS Computational Biology}, 7(3):e1002008, 2011.

\bibitem{giunta2022local}
Valeria Giunta, Thomas Hillen, Mark Lewis, and Jonathan~R Potts.
\newblock Local and global existence for nonlocal multispecies
  advection-diffusion models.
\newblock {\em SIAM Journal on Applied Dynamical Systems}, 21(3):1686--1708,
  2022.

\bibitem{GHLP22}
Valeria Giunta, Thomas Hillen, Mark~A Lewis, and Jonathan~R Potts.
\newblock Detecting minimum energy states and multi-stability in nonlocal
  advection--diffusion models for interacting species.
\newblock {\em Journal of Mathematical Biology}, 85(5):1--44, 2022.

\bibitem{giunta2021pattern}
Valeria Giunta, Maria~Carmela Lombardo, and Marco Sammartino.
\newblock Pattern formation and transition to chaos in a chemotaxis model of
  acute inflammation.
\newblock {\em SIAM Journal on Applied Dynamical Systems}, 20(4):1844--1881,
  2021.

\bibitem{goldstone2009collective}
Robert~L Goldstone and Todd~M Gureckis.
\newblock Collective behavior.
\newblock {\em Topics in Cognitive Science}, 1(3):412--438, 2009.

\bibitem{usersguide}
Thomas Hillen and Kevin~J Painter.
\newblock A user’s guide to pde models for chemotaxis.
\newblock {\em Journal of {M}athematical {B}iology}, 58(1-2):183, 2009.

\bibitem{hoyle06}
Rebecca Hoyle.
\newblock {\em Pattern formation: an introduction to methods}.
\newblock Cambridge University Press, 2006.

\bibitem{jungel2022nonlocal}
Ansgar J{\"u}ngel, Stefan Portisch, and Antoine Zurek.
\newblock Nonlocal cross-diffusion systems for multi-species populations and
  networks.
\newblock {\em Nonlinear Analysis}, 219:112800, 2022.

\bibitem{MC00}
PC~Matthews and Stephen~M Cox.
\newblock Pattern formation with a conservation law.
\newblock {\em Nonlinearity}, 13(4):1293, 2000.

\bibitem{Murray01}
James~D Murray.
\newblock {\em Mathematical biology II: Spatial models and biomedical
  applications}, volume~3.
\newblock Springer New York, 2001.

\bibitem{potts2014animal}
Jonathan~R Potts and Mark~A Lewis.
\newblock How do animal territories form and change? lessons from 20 years of
  mechanistic modelling.
\newblock {\em Proceedings of the Royal Society B: Biological Sciences},
  281(1784):20140231, 2014.

\bibitem{potts2016territorial}
Jonathan~R Potts and Mark~A Lewis.
\newblock Territorial pattern formation in the absence of an attractive
  potential.
\newblock {\em Journal of Mathematical Biology}, 72:25--46, 2016.

\bibitem{PL19}
Jonathan~R Potts and Mark~A Lewis.
\newblock Spatial memory and taxis-driven pattern formation in model
  ecosystems.
\newblock {\em Bulletin of mathematical biology}, 81(7):2725--2747, 2019.

\bibitem{puckett2014searching}
James~G Puckett, Douglas~H Kelley, and Nicholas~T Ouellette.
\newblock Searching for effective forces in laboratory insect swarms.
\newblock {\em Scientific Reports}, 4(1):1--5, 2014.

\bibitem{rodriguez2020steady}
Nancy Rodr{\'\i}guez and Yi~Hu.
\newblock On the steady-states of a two-species non-local cross-diffusion
  model.
\newblock {\em Journal of Applied Analysis}, 26(1):1--19, 2020.

\bibitem{soresina2022hopf}
Cinzia Soresina.
\newblock Hopf bifurcations in the full skt model and where to find them.
\newblock {\em arXiv preprint arXiv:2202.04168}, 2022.

\bibitem{tao2016dynamic}
Yun Tao, Luca B{\"o}rger, and Alan Hastings.
\newblock Dynamic range size analysis of territorial animals: An optimality
  approach.
\newblock {\em The American Naturalist}, 188(4):460--474, 2016.

\bibitem{TBL06}
Chad~M Topaz, Andrea~L Bertozzi, and Mark~A Lewis.
\newblock A nonlocal continuum model for biological aggregation.
\newblock {\em Bulletin of Mathematical Biology}, 68:1601--1623, 2006.

\bibitem{TB90}
Laurette~S Tuckerman and Dwight Barkley.
\newblock Bifurcation analysis of the eckhaus instability.
\newblock {\em Physica D: Nonlinear Phenomena}, 46(1):57--86, 1990.

\bibitem{Turing1952}
Alan~M Turing.
\newblock The chemical basis of morphogenesis.
\newblock {\em Philosophical Transactions of the Royal Society of London.
  Series B, Biological Sciences}, 237(641):37--72, 1952.

\bibitem{uecker2021numerical}
Hannes Uecker.
\newblock {\em Numerical continuation and bifurcation in Nonlinear PDEs}.
\newblock SIAM, 2021.

\bibitem{verschueren2017model}
Nicolas Verschueren and Alan Champneys.
\newblock A model for cell polarization without mass conservation.
\newblock {\em SIAM Journal on Applied Dynamical Systems}, 16(4):1797--1830,
  2017.

\end{thebibliography}



\end{document}